# *An Existence Theory for Small-Amplitude Doubly Periodic Water Waves with Vorticity*

E. Lokharu, D. S. Seth & E. Wahlén

*Communicated by* N. Masmoudi

**Abstract**

We prove the existence of three-dimensional steady gravity-capillary waves with vorticity on water of finite depth. The waves are periodic with respect to a given two-dimensional lattice and the relative velocity field is a Beltrami field, meaning that the vorticity is collinear to the velocity. The existence theory is based on multi-parameter bifurcation theory.

## 1. Introduction

### *1.1. Statement of the Problem*

This paper is concerned with three-dimensional steady water waves driven by gravity and surface tension. An inviscid fluid of constant unit density occupies the domain

$$\Omega^\eta = \{(\boldsymbol{x}', z) \in \mathbb{R}^2 \times \mathbb{R} : -d < z < \eta(\boldsymbol{x}')\}$$

for some function $\eta\colon \mathbb{R}^2 \to \mathbb{R}$ and constant $d > 0$, where $\boldsymbol{x}' = (x, y)$. Let $\boldsymbol{u}\colon \overline{\Omega^\eta} \to \mathbb{R}^3$ be the (relative) velocity field and $p\colon \overline{\Omega^\eta} \to \mathbb{R}$ the pressure. In a frame of reference moving with the wave, the fluid motion is governed by the stationary Euler equations

$$\begin{aligned} (\boldsymbol{u} \cdot \nabla)\boldsymbol{u} &= -\nabla p - g e_3 && \text{in } \Omega^\eta, \\ \nabla \cdot \boldsymbol{u} &= 0 && \text{in } \Omega^\eta, \end{aligned}$$

with a kinematic boundary condition on the top and bottom boundaries:

$$\boldsymbol{u} \cdot \boldsymbol{n} = 0 \qquad \text{on } \partial\Omega^\eta,$$

and a dynamic boundary condition on the free surface:

$$p = p_{atm} - 2\sigma K_M \quad \text{on } z = \eta(\boldsymbol{x}').$$

Here $\boldsymbol{e}_3 = (0, 0, 1)$, $K_M$ is the mean curvature of the free surface, given by

$$2K_M = \nabla \cdot \left( \frac{\nabla \eta}{\sqrt{1 + |\nabla \eta|^2}} \right),$$

while $\sigma > 0$ is the coefficient of surface tension and $p_{atm}$ the constant atmospheric pressure. Supposing that the wave is moving with constant velocity $\boldsymbol{\nu} = (\nu_1, \nu_2)$ in the original stationary frame of reference, then in this frame the travelling wave is given by $z = \eta(\boldsymbol{x}' - \boldsymbol{\nu}t)$ and the corresponding velocity field is given by $\boldsymbol{v} = \boldsymbol{u}(\boldsymbol{x}' - \boldsymbol{\nu}t, z) + (\boldsymbol{\nu}, 0)$. Note that $\boldsymbol{v}$ satisfies $\nabla \times \boldsymbol{v} = \alpha \boldsymbol{v} - \alpha(\boldsymbol{\nu}, 0)$ and is therefore not a Beltrami field in general.

Almost all previous studies of three-dimensional steady water waves have been restricted to the irrotational setting, where $\nabla \times \boldsymbol{u} = 0$. It is desirable to relax this condition in order to model interactions of surface waves with non-uniform currents. In the present paper we consider the special case when the velocity and vorticity fields are collinear, that is,

$$\nabla \times \boldsymbol{u} = \alpha \boldsymbol{u} \quad \text{in } \Omega^\eta$$

for some constant $\alpha$. In other words, we assume that $\boldsymbol{u}$ is a (strong) Beltrami field. Such fields are well-known in solar and plasma physics (see e.g. BOULMEZAOUD, MADAY & AMARI [5], FREIDBERG [16] and PRIEST [30]) and are also called (linear) force-free fields. The adjectives 'strong' and 'linear' refer to the fact that $\alpha$ is assumed to be constant. The more complicated case when $\alpha$ is variable has been investigated by several authors (see e.g. BOULMEZAOUD & AMARI [4], ENCISO & PERALTA-SALAS [15] and KAISER, NEUDERT & VON WAHL [26]), but will not be considered herein. Any divergence-free Beltrami field generates a solution to the stationary Euler equations with pressure

$$p = C - \frac{|\boldsymbol{u}|^2}{2} - gz.$$

The governing equations are thus replaced by

$$\nabla \times \boldsymbol{u} = \alpha \boldsymbol{u} \quad \text{in } \Omega^\eta, \tag{1.1a}$$

$$\nabla \cdot \boldsymbol{u} = 0 \quad \text{in } \Omega^\eta, \tag{1.1b}$$

$$\boldsymbol{u} \cdot \boldsymbol{n} = 0 \quad \text{on } \partial\Omega^\eta, \tag{1.1c}$$

$$\frac{1}{2}|\boldsymbol{u}|^2 + g\eta - \sigma \nabla \cdot \left( \frac{\nabla \eta}{\sqrt{1 + |\nabla \eta|^2}} \right) = Q \quad \text{on } z = \eta, \tag{1.1d}$$

where $Q$ is the Bernoulli constant. Condition (1.1b) is actually redundant for $\alpha \neq 0$, but we retain it since we want to allow that $\alpha = 0$. For a given $\eta$ there can be more than one solution to the above equations and we will therefore later append integral conditions in order to enforce uniqueness.

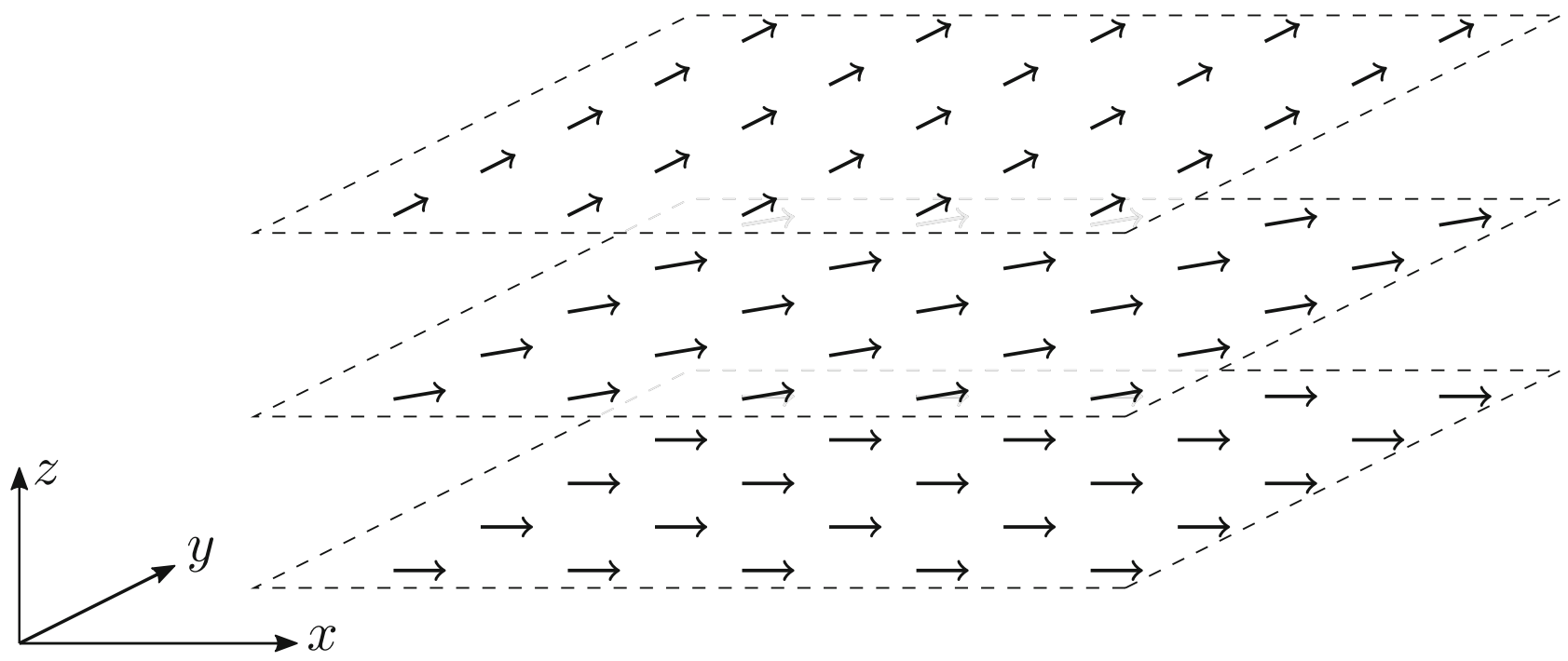

**Fig. 1.** A laminar flow in different horizontal sections of the fluid domain

The choice of Beltrami flows is mainly motivated by mathematical considerations, since it gives rise to an elliptic free boundary problem. From a physical point of view, the choice is quite specific and it would be desirable to treat more general flows. One interesting feature of Beltrami flows is that they include laminar flows whose direction varies with depth (see Section 1.2.1 and Figure 1). This could potentially be of interest when considering a wind-induced surface current interacting with a subsurface current in a different direction.

## *1.2. Special Solutions*

**1.2.1. Laminar Flows** Let us consider a fluid domain with a flat boundary, that is $\eta \equiv 0$:

$$\Omega^0 = \{(\boldsymbol{x}', z) \in \mathbb{R}^3 : -d < z < 0,\ \boldsymbol{x}' \in \mathbb{R}^2\}.$$

In this case we find a two-parameter family of ‘trivial’ solutions given by laminar flows

$$\boldsymbol{U}[c_1, c_2] = c_1\boldsymbol{U}^{(1)} + c_2\boldsymbol{U}^{(2)}, \quad c_1, c_2 \in \mathbb{R},$$

where

$$\boldsymbol{U}^{(1)} = (\cos(\alpha z), -\sin(\alpha z), 0), \quad \boldsymbol{U}^{(2)} = (\sin(\alpha z), \cos(\alpha z), 0),$$

and the corresponding Bernoulli constant $Q$ in (1.1d) is given by

$$Q(c_1, c_2) = \frac{1}{2}[c_1^2 + c_2^2].$$

The laminar flow $\boldsymbol{U}$ is constant in every horizontal section of the fluid domain but the direction of the flow depends on the vertical coordinate (see Figure 1). The constants $c_1, c_2$ will be used later as bifurcation parameters.

**1.2.2. Two-and-a-Half-Dimensional Waves** There is a connection between problem (1.1) and the two-dimensional steady water wave problem with affine linear vorticity function. Indeed, let $\eta(\bar{x})$ be the surface and $\psi(\bar{x}, z)$ be the stream function for a two-dimensional wave with vorticity function $\alpha^2\psi + \alpha\beta$ traveling in the $\bar{\boldsymbol{e}}$-direction, where $\alpha$ and $\beta$ are real constants, $\bar{\boldsymbol{e}}$ is a horizontal unit vector and $\bar{x} = \boldsymbol{x} \cdot \bar{\boldsymbol{e}}$. Then

$$\begin{aligned} &\partial_{\bar{x}}^2\psi + \partial_z^2\psi + \alpha^2\psi + \alpha\beta = 0 \ \text{ in } \ \Omega_{2D}^{\eta}, \\ &\psi(\bar{x}, -d) = m_1, \quad \psi(\bar{x}, \eta(\bar{x})) = m_2, \quad \bar{x} \in \mathbb{R}, \\ &\frac{1}{2}|\nabla\psi|^2 + gz - \sigma\partial_{\bar{x}}\left(\frac{\partial_{\bar{x}}\eta}{\sqrt{1 + (\partial_{\bar{x}}\eta)^2}}\right) = Q_0, \quad z = \eta(\bar{x}), \end{aligned} \tag{1.2}$$

where $m_1, m_2, Q_0$ are constants, while

$$\Omega_{2D}^{\eta} = \{(\bar{x}, z) \in \mathbb{R}^2 : \ -d < z < \eta(\bar{x})\}$$

is the two-dimensional fluid region. The corresponding velocity field is given by

$$\boldsymbol{u}(\bar{x}, z) = -\psi_z(\bar{x}, z)\bar{\boldsymbol{e}} + \psi_{\bar{x}}(\bar{x}, z)\boldsymbol{e}_3.$$

We can turn this into a solution of the stationary Euler equations in the three-dimensional domain by letting $\eta$ and $\boldsymbol{u}$ equal the two-dimensional solution for every $\bar{y} = \boldsymbol{x} \cdot \boldsymbol{e}_\perp$, where $\boldsymbol{e}_\perp = \boldsymbol{e}_3 \times \bar{\boldsymbol{e}}$. However, this solution is clearly still two-dimensional in the sense that it is independent of $\bar{y}$ and the velocity vector is collinear with the direction of propagation, and hence there is no fluid motion in the perpendicular horizontal direction $\boldsymbol{e}_\perp$. On the other hand, we can put

$$\boldsymbol{u}(\bar{x}, \bar{y}, z) = -\psi_z(\bar{x}, z)\bar{\boldsymbol{e}} + (\alpha\psi(\bar{x}, z) + \beta)\boldsymbol{e}_\perp + \psi_{\bar{x}}(\bar{x}, z)\boldsymbol{e}_3. \tag{1.3}$$

One verifies that $\boldsymbol{u}$ solves (1.1) in

$$\Omega^{\eta} := \{\bar{x}\bar{\boldsymbol{e}} + \bar{y}\boldsymbol{e}_\perp + z\boldsymbol{e}_3 : -d < z < \eta(\bar{x}), \quad \bar{x}, \bar{y} \in \mathbb{R}\}.$$

The flow generated by $\boldsymbol{u}$ is called 2½-dimensional, since $\boldsymbol{u}$ only depends on the two variables $\bar{x}$ and $z$ but has a non-zero $\bar{y}$-component; see MAJDA & BERTOZZI [28, Sect. 2.3]. Note that every laminar solution $\boldsymbol{U}[c_1, c_2]$ can be written in the form (1.3) for some stream function $\Psi(z)$. Note also that one can get rid of the constant $\beta$ when $\alpha \neq 0$ by introducing the new stream function $\psi + \alpha^{-1}\beta$.

Conversely, any 2½-dimensional Beltrami flow arises from a solution to the two-dimensional steady water wave problem with affine linear vorticity function. Indeed, assume that we have a solution $(\boldsymbol{u}, \eta)$ to (1.1) depending only on one horizontal variable $\bar{x}$. Then $(\bar{u}, u_3)$ is divergence free with respect to the variables $(\bar{x}, z)$, where $\bar{u} = \boldsymbol{u} \cdot \bar{\boldsymbol{e}}$, and hence there exists a stream function $\psi(\bar{x}, z)$ such that $u_3 = \psi_{\bar{x}}$ and $\bar{u} = -\psi_z$. Now, equation (1.1a) gives $u_\perp = \alpha\psi + \beta$ for some constant $\beta$, as well as

$$\psi_{\bar{x}\bar{x}} + \psi_{zz} + \alpha^2\psi + \alpha\beta = 0.$$

On the other hand, $\boldsymbol{u}$ is subject to (1.1c), which implies that $\psi$ is constant on the boundaries. Using this fact we recover the Bernoulli equation for $\psi$ from (1.1d).

Problem (1.2) with zero surface tension has for example been studied by AASEN & VARHOLM [1], EHRNSTRÖM, ESCHER & WAHLÉN [13], EHRNSTRÖM, ESCHER & VILLARI [12] and EHRNSTRÖM & WAHLÉN [14]. The gravity-capillary problem has been considered for a general class of vorticity functions (including affine) but restricted to flows without stagnation points by several authors; see e.g. WAHLÉN [33] and WALSH [35,36]. These two-dimensional existence results immediately yield the existence of 2½-dimensional waves on Beltrami flows. In this paper we will instead take the opposite approach. As a part of the analysis we will directly obtain the existence of 2½-dimensional waves, which generate solutions of problem (1.2) by the above correspondence; see Remark 4.7.

### *1.3. Previous Results*

The theory of three-dimensional steady waves with vorticity is a relatively new subject of studies. In contrast to the two-dimensional case (see e.g. CONSTANTIN [9]) one cannot, in general, reformulate the problem as an elliptic free boundary problem. At the moment there are no existence results, save for some explicit Gerstner-type solutions for edge waves along a sloping beach and equatorially trapped waves; see CONSTANTIN [8,10] and HENRY [23] and references therein. WAHLÉN [34] showed that the assumption of constant vorticity prevents the existence of genuinely three-dimensional traveling gravity waves on water of finite depth. A variational principle for doubly periodic waves whose relative velocity is given by a Beltrami vector field was obtained by LOKHARU & WAHLÉN [27].

The irrotational theory is on the other hand much more developed. The first existence proofs for doubly periodic, irrotational, gravity-capillary waves in the 'non-resonant' case are due to REEDER & SHINBROT [31] and SUN [32]. These papers consider periodic lattices for which the fundamental domain is a 'symmetric diamond'. The resonant case was investigated by CRAIG & NICHOLLS [11] using a combination of topological and variational methods. They proved the existence of small-amplitude periodic waves for an arbitrary fundamental domain. A different approach known as spatial dynamics was developed by GROVES & MIELKE [20]. The idea is to choose one spatial variable for the role of time and think of the problem as a Hamiltonian system with an infinite dimensional phase space. Using this approach, GROVES & MIELKE constructed symmetric doubly periodic waves. The asymmetric case was later investigated by GROVES & HARAGUS [18]; see also NILSSON [29]. One of the strengths of spatial dynamics is that is not restricted to the doubly periodic setting. It can also be used to construct waves which e.g. are solitary in one direction and periodic or quasi-periodic in another; see GROVES [17] for a survey of different results. One type of solutions which have so far elluded the spatial dynamics method is fully localised solitary waves, that is, solutions which decay in all horizontal directions. Such solutions have however been constructed using variational methods; see BUFFONI, GROVES, SUN & WAHLÉN [6], BUFFONI, GROVES & WAHLÉN [7] and GROVES & SUN [21]. Finally, note that the doubly periodic problem with zero surface tension is considerably harder since one has to deal with small divisors. Nevertheless, IOOSS & PLOTNIKOV [24,25] proved existence results for symmetric and asymmetric waves using Nash-Moser techniques. It might be possible to deal with the zero surface tension version of problem (1.1) in a similar way.

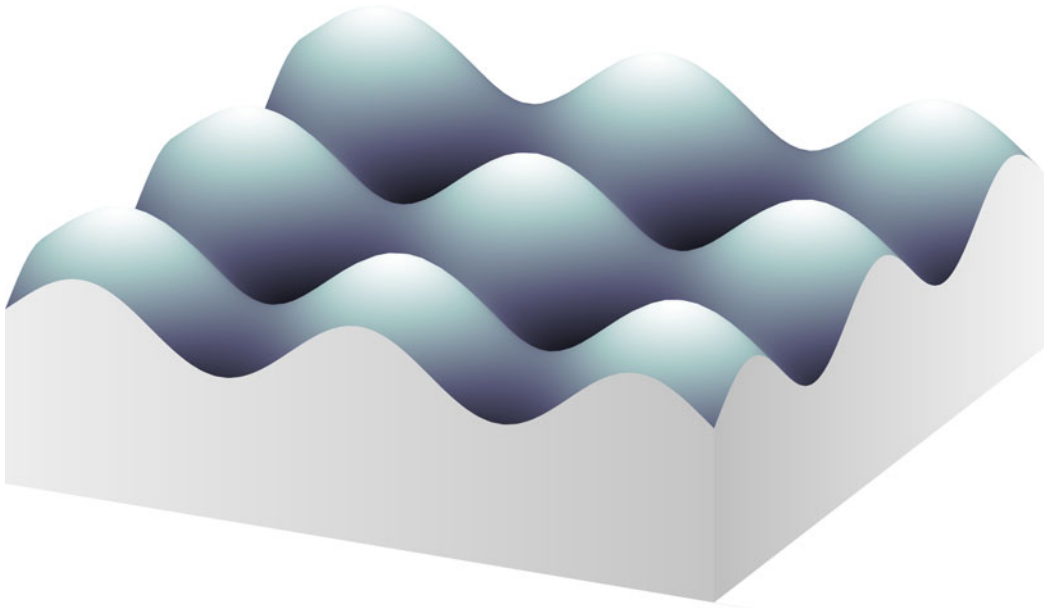

**Fig. 2.** A sketch of a doubly periodic wave

### 1.4. The Present Contribution

The main contribution of our paper is an existence result for small-amplitude doubly periodic solutions of problem (1.1). Given two linearly independent vectors $\boldsymbol{\lambda}_1, \boldsymbol{\lambda}_2 \in \mathbb{R}^2$ we define the two-dimensional lattice

$$\Lambda = \{\boldsymbol{\lambda} = m_1\boldsymbol{\lambda}_1 + m_2\boldsymbol{\lambda}_2 : m_1, m_2 \in \mathbb{Z}\}.$$

We assume that

$$\eta(\boldsymbol{x}' + \boldsymbol{\lambda}) = \eta(\boldsymbol{x}') \tag{1.4a}$$

and

$$\boldsymbol{u}(\boldsymbol{x}' + \boldsymbol{\lambda}, z) = \boldsymbol{u}(\boldsymbol{x}', z) \tag{1.4b}$$

for $\boldsymbol{\lambda} \in \Lambda$, so that the fluid domain $\Omega^\eta$ and velocity field are periodic with respect to the lattice $\Lambda$ (see Figure 2). In addition, we impose the symmetry conditions

$$\eta(-\boldsymbol{x}') = \eta(\boldsymbol{x}'), \tag{1.5a}$$

$$\boldsymbol{u}(-\boldsymbol{x}', z) = (u_1(\boldsymbol{x}', z), u_2(\boldsymbol{x}', z), -u_3(\boldsymbol{x}', z)). \tag{1.5b}$$

For later use it is convenient to define

$$B_{lj} = \{(a_1 + l)\boldsymbol{\lambda}_1 + (a_2 + j)\boldsymbol{\lambda}_2 : a_1, a_2 \in (0, 1)\}, \quad l, j \in \mathbb{Z},$$

and

$$\Omega^\eta_{lj} = \{(\boldsymbol{x}', z) : -d < z < \eta(\boldsymbol{x}'),\ \boldsymbol{x}' \in B_{lj}\}, \quad l, j \in \mathbb{Z},$$

which splits the domain $\Omega^\eta$ into simple periodic cells. We denote the top and bottom boundaries of $\Omega^\eta_{lj}$ by $\partial\Omega^{\eta,s}_{lj}$ and $\partial\Omega^{\eta,b}_{lj}$ respectively (see Figure 3).

We study solutions bifurcating from laminar flows $\boldsymbol{U}[c_1, c_2]$, where $c_1$ and $c_2$ act as bifurcation parameters and therefore vary along the family of nontrivial solutions that we find. We look for solutions satisfying (1.1d) with the same constant $Q = Q(c_1, c_2)$ as the underlying laminar flow $\boldsymbol{U}[c_1, c_2]$. Therefore, the Bernoulli constant $Q$ will vary along the bifurcation set. For purposes of uniqueness we also impose integral conditions relating the total (relative) horizontal momentum in the $x$ and $y$ directions to that of $\boldsymbol{U}[c_1, c_2]$. This results in the system

$$\nabla \times \boldsymbol{u} = \alpha \boldsymbol{u} \quad \text{in } \Omega^\eta, \tag{1.6a}$$

$$\nabla \cdot \boldsymbol{u} = 0 \quad \text{in } \Omega^\eta, \tag{1.6b}$$

$$\boldsymbol{u} \cdot \boldsymbol{n} = 0 \quad \text{on } \partial\Omega^\eta, \tag{1.6c}$$

$$\int_{\Omega_{00}^\eta} u_j \, \mathrm{d}V = \int_{\Omega_{00}^\eta} U_j[c_1, c_2] \, \mathrm{d}V \quad j = 1, 2, \tag{1.6d}$$

$$\frac{1}{2}|\boldsymbol{u}|^2 + g\eta - \sigma \nabla \cdot \left( \frac{\nabla \eta}{\sqrt{1 + |\nabla \eta|^2}} \right) = Q(c_1, c_2) \quad \text{on } z = \eta. \tag{1.6e}$$

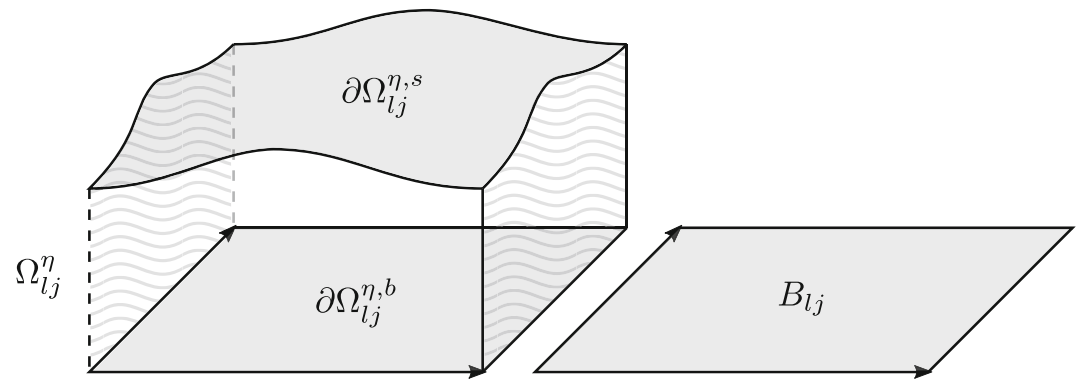

**Fig. 3.** A three-dimensional periodic cell of the domain (left) and the corresponding two-dimensional periodic cell (right)

In Section 2 we introduce a suitable functional-analytic framework. Since we are dealing with a free boundary problem it is convenient to transform the problem to a fixed domain. After a sequence of changes of variables we derive an equivalent version, problem (2.5), which is amenable to further analysis. We also show that it is possible to reduce the governing equations to a single nonlinear pseudodifferential equation for the free surface (cf. Theorem 2.1 and equation (2.6)). In Section 3 we study the linearisation of the problem and identify bifurcation points $(c_1^\star, c_2^\star)$ at which its solution space is two-dimensional. In doing so, we derive a dispersion relation $\rho(\boldsymbol{c}, \boldsymbol{k}) = 0$ (see equation (3.3)) which shows how the parameters $\boldsymbol{c} = (c_1, c_2)$ are related to the wave vector $\boldsymbol{k}$ of a solution to the linearised problem. The bifurcation points are those for which $\rho(\boldsymbol{c}^\star, \boldsymbol{k}_1) = \rho(\boldsymbol{c}^\star, \boldsymbol{k}_2) = 0$ for two linearly independent wave vectors $\boldsymbol{k}_1, \boldsymbol{k}_2$. The conclusions can be found in Propositions 3.1 and 3.3. In Section 4 we finally give a precise formulation and proof of the main result, Theorem 4.1. While the main interest of this paper lies in the case $\alpha \neq 0$, the existence result also covers the case $\alpha = 0$ and therefore yields another existence proof for doubly periodic irrotational gravity-capillary waves. In the irrotational case there is an additional symmetry which allows one to treat the case of symmetric fundamental domains using classical one-dimensional bifurcation theory. The lack of this symmetry is one of the reasons for using a multi-parameter bifurcation approach. We have formulated the main result in terms of the relative velocity field. It is worth keeping in mind that in the original stationary frame the solution is a small perturbation of the laminar flow $\boldsymbol{U} + (\boldsymbol{v}, 0)$, where $\boldsymbol{v}$ is the velocity vector.

## 2. Functional-Analytic Framework

### 2.1. Function Spaces and Notation

Suppose that $\Omega$ is an open subset of Euclidean space and let $C^{k,\gamma}(\overline{\Omega})$, with $k \in \mathbb{N}_0 := \{0, 1, 2, \ldots\}$ and $\gamma \in (0, 1)$, denote the class of functions $u \colon \overline{\Omega} \to \mathbb{R}$ whose partial derivatives up to order $k$ are bounded and uniformly $\gamma$-Hölder continuous in $\overline{\Omega}$. This is a Banach space when equipped with the norm

$$\|u\|_{C^{k,\gamma}(\overline{\Omega})} := \max_{|\beta|\le k} \sup_{\boldsymbol{x}\in\overline{\Omega}} |\partial^\beta u(\boldsymbol{x})| + \max_{|\beta|=k} [\partial^\beta u]_{0,\gamma},$$

with

$$[v]_{0,\gamma} := \sup_{\boldsymbol{x}\neq\boldsymbol{y}\in\overline{\Omega}} \frac{|v(\boldsymbol{x}) - v(\boldsymbol{y})|}{|\boldsymbol{x} - \boldsymbol{y}|^\gamma}$$

and $\beta$ denoting multiindices. We will consider surface profiles in the space $C^{k,\gamma}_{per,e}(\mathbb{R}^2)$, consisting of $\eta \in C^{k,\gamma}(\mathbb{R}^2)$ which satisfy the periodicity condition (1.4a) and the evenness condition (1.5a), and velocity fields in the space $(C^{k,\gamma}_{per,e}(\overline{\Omega^\eta}))^2 \times C^{k,\gamma}_{per,o}(\overline{\Omega^\eta})$ consisting of $\boldsymbol{u} \in (C^{k,\gamma}(\overline{\Omega^\eta}))^3$ which satisfy (1.4b) and (1.5b).

We let

$$\Lambda' := \{\boldsymbol{k} = n_1\boldsymbol{k}_1 + n_2\boldsymbol{k}_2 : n_1, n_2 \in \mathbb{Z},\ \boldsymbol{k}_i \cdot \boldsymbol{\lambda}_j = 2\pi\delta_{ij},\ i, j = 1, 2\}$$

denote the lattice dual to $\Lambda$. Note that any $\eta \in C^{k,\gamma}_{per,e}(\mathbb{R}^2)$ can be expanded in a Fourier series

$$\eta(\boldsymbol{x}') = \sum_{\boldsymbol{k}\in\Lambda'} \hat{\eta}^{(\boldsymbol{k})} e^{i\boldsymbol{k}\cdot\boldsymbol{x}'},$$

with Fourier coefficients

$$\hat{\eta}^{(\boldsymbol{k})} = \frac{1}{|B_{00}|} \int_{B_{00}} \eta(\boldsymbol{x}') e^{i\boldsymbol{k}\cdot\boldsymbol{x}'} \,\mathrm{d}x\,\mathrm{d}y$$

satisfying $\hat{\eta}^{(\boldsymbol{k})} = \hat{\eta}^{(-\boldsymbol{k})} \in \mathbb{R}$. If on the other hand $\eta \in C^{k,\gamma}_{per,o}(\mathbb{R}^2)$, then $\hat{\eta}^{(\boldsymbol{k})}$ is purely imaginary with $\hat{\eta}^{(\boldsymbol{k})} = -\hat{\eta}^{(-\boldsymbol{k})}$. Functions in $C^{k,\gamma}_{per,e}(\overline{\Omega^0})$ or $C^{k,\gamma}_{per,o}(\overline{\Omega^0})$ have a similar expansion, with Fourier coefficients depending on $z$. Note that all of the analysis can also be done in Sobolev spaces.

If $X$ and $Y$ are normed vector spaces and $G \colon X \to Y$ is a Fréchet differentiable mapping, we will denote its Fréchet derivative at $x \in X$ by $DG[x]$.

*2.2. Flattening Transformation*

Since (1.6) is a free boundary problem, it is convenient to perform a change of variables which fixes the domain. Under the ‘flattening’ transformation

$$\boldsymbol{F}\colon (\dot{x}, \dot{y}, \dot{z}) \mapsto (x, y, z) = \left(\dot{x}, \dot{y}, \dot{z} + \eta(\dot{x}, \dot{y}) \left(\frac{\dot{z}}{d} + 1\right)\right),$$

the fluid domain $\Omega^\eta$ becomes the image of the ‘flat’ domain $\Omega^0$. In the new variables $\dot{\boldsymbol{x}} = (\dot{x}, \dot{y}, \dot{z})$ we introduce the position dependent basis $\boldsymbol{f}_1(\dot{\boldsymbol{x}}) = \frac{\partial \boldsymbol{F}(\dot{\boldsymbol{x}})}{\partial \dot{x}}$, $\boldsymbol{f}_2(\dot{\boldsymbol{x}}) = \frac{\partial \boldsymbol{F}(\dot{\boldsymbol{x}})}{\partial \dot{y}}$, $\boldsymbol{f}_3(\dot{\boldsymbol{x}}) = \frac{\partial \boldsymbol{F}(\dot{\boldsymbol{x}})}{\partial \dot{z}}$ given explicitly by

$$\boldsymbol{f}_1 = \begin{pmatrix} 1 \\ 0 \\ \eta_x \frac{\dot{z}+d}{d} \end{pmatrix}, \quad \boldsymbol{f}_2 = \begin{pmatrix} 0 \\ 1 \\ \eta_y \frac{\dot{z}+d}{d} \end{pmatrix}, \quad \boldsymbol{f}_3 = \begin{pmatrix} 0 \\ 0 \\ \frac{\eta+d}{d} \end{pmatrix}.$$

Note that the vectors $\boldsymbol{f}_j(\dot{\boldsymbol{x}})$ are tangent to the coordinate curves. In the new variables, the vector field $\boldsymbol{u}$ is given by

$$\dot{\boldsymbol{u}}(\dot{\boldsymbol{x}}) = (\dot{u}_1(\dot{\boldsymbol{x}}), \dot{u}_2(\dot{\boldsymbol{x}}), \dot{u}_3(\dot{\boldsymbol{x}})),$$

where the coordinate functions are determined by

$$\begin{aligned} &\dot{u}_1(\dot{\boldsymbol{x}})\boldsymbol{f}_1(\dot{\boldsymbol{x}}) + \dot{u}_2(\dot{\boldsymbol{x}})\boldsymbol{f}_2(\dot{\boldsymbol{x}}) + \dot{u}_3(\dot{\boldsymbol{x}})\boldsymbol{f}_3(\dot{\boldsymbol{x}}) \\ &\quad = J(\boldsymbol{F}^{-1}(\boldsymbol{x}))[u_1(\boldsymbol{x})\boldsymbol{e}_1 + u_2(\boldsymbol{x})\boldsymbol{e}_2 + u_3(\boldsymbol{x})\boldsymbol{e}_3], \end{aligned}$$

in which

$$J(\dot{\boldsymbol{x}}) = \det \boldsymbol{F}'(\dot{\boldsymbol{x}}) = \frac{\eta(\dot{\boldsymbol{x}}) + d}{d}.$$

Note that under this change of variables the space $(C^{k,\gamma}_{per,e}(\overline{\Omega^\eta}))^2 \times C^{k,\gamma}_{per,o}(\overline{\Omega^\eta})$ is mapped bijectively to $(C^{k,\gamma}_{per,e}(\overline{\Omega^0}))^2 \times C^{k,\gamma}_{per,o}(\overline{\Omega^0})$ if $\eta \in C^{k+1,\gamma}_{per,e}(\mathbb{R}^2)$ with $\min \eta > -d$. The divergence and curl take the following forms in the new coordinates:

$$\nabla \cdot \boldsymbol{u} = [J(\dot{\boldsymbol{x}})]^{-1} \nabla_{\dot{\boldsymbol{x}}} \cdot \dot{\boldsymbol{u}} \tag{2.1}$$

and

$$\nabla \times \boldsymbol{u} = \frac{1}{J(\dot{\boldsymbol{x}})} \det \begin{pmatrix} \boldsymbol{f}_1 & \boldsymbol{f}_2 & \boldsymbol{f}_3 \\ \partial_{\dot{x}} & \partial_{\dot{y}} & \partial_{\dot{z}} \\ \boldsymbol{f}_1 \cdot \boldsymbol{u}(\boldsymbol{F}(\dot{\boldsymbol{x}})) & \boldsymbol{f}_2 \cdot \boldsymbol{u}(\boldsymbol{F}(\dot{\boldsymbol{x}})) & \boldsymbol{f}_3 \cdot \boldsymbol{u}(\boldsymbol{F}(\dot{\boldsymbol{x}})) \end{pmatrix}. \tag{2.2}$$

Thus, in view of (2.1) and (2.2) problem (1.6) transforms into

$$\nabla \times \dot{\boldsymbol{u}} - \alpha \dot{\boldsymbol{u}} = \nabla \times \boldsymbol{N}(\dot{\boldsymbol{u}}, \eta) \qquad \text{in } \Omega^0, \tag{2.3a}$$

$$\nabla \cdot \dot{\boldsymbol{u}} = 0 \qquad \text{in } \Omega^0, \tag{2.3b}$$

$$\dot{u}_3 = 0 \qquad \text{on } \partial\Omega^0, \tag{2.3c}$$

$$\int_{\Omega^0_{00}} \dot{u}_j \, \mathrm{d}V = \int_{\Omega^\eta_{00}} U_j[c_1, c_2] \, \mathrm{d}V \qquad j = 1, 2, \tag{2.3d}$$

$$\frac{1}{2}B(\dot{\boldsymbol{u}}, \eta) + g\eta - \sigma\nabla \cdot \left( \frac{\nabla\eta}{\sqrt{1+|\nabla\eta|^2}} \right) = Q(c_1, c_2) \qquad \text{on } \dot{z} = 0, \tag{2.3e}$$

where the divergence and curl are now with respect to the dotted variables and the nonlinearities $B$ and $\boldsymbol{N} = (N_1, N_2, N_3)$ are given by

$$N_j(\dot{\boldsymbol{u}}, \eta) = \dot{u}_j - \sum_{l=1}^{3} \frac{\boldsymbol{f}_j \cdot \boldsymbol{f}_l}{J} \dot{u}_l, \quad j = 1, 2, 3,$$

$$B(\dot{\boldsymbol{u}}, \eta) = \frac{1}{J^2}\left[ \dot{u}_1^2 + \dot{u}_2^2 + \left( \eta_x \frac{\dot{z}+d}{d}\dot{u}_1 + \eta_y \frac{\dot{z}+d}{d}\dot{u}_2 + \frac{\eta+d}{d}\dot{u}_3 \right)^2 \right].$$

Note that $\boldsymbol{N}$ is linear in $\dot{\boldsymbol{u}}$ and that $\boldsymbol{N}(\dot{\boldsymbol{u}}, 0) = 0$. From now on we drop the dots on $x$, $y$ and $z$ to simplify the notation.

Since we are interested in solutions that are close to a laminar flow $\boldsymbol{U}[c_1, c_2]$, we write

$$\dot{\boldsymbol{u}} = \boldsymbol{U} + \tilde{\boldsymbol{v}}.$$

Thus, equations (2.3) transform into

$$\nabla \times \tilde{\boldsymbol{v}} - \alpha\tilde{\boldsymbol{v}} = \nabla \times \tilde{\boldsymbol{L}}(\eta) + \nabla \times \tilde{\boldsymbol{N}}(\tilde{\boldsymbol{v}}, \eta) \qquad \text{in } \Omega^0, \tag{2.4a}$$

$$\nabla \cdot \tilde{\boldsymbol{v}} = 0 \qquad \text{in } \Omega^0, \tag{2.4b}$$

$$v_3 = 0 \qquad \text{on } \partial\Omega^0, \tag{2.4c}$$

$$\int_{\Omega_{00}^0} \tilde{v}_j \,\mathrm{d}V = \int_{\Omega_{00}^\eta} U_j \,\mathrm{d}V - \int_{\Omega_{00}^0} U_j \,\mathrm{d}V \qquad j = 1, 2, \tag{2.4d}$$

$$-[c_1^2 + c_2^2]\frac{\eta}{d} + c_1\tilde{v}_1 + c_2\tilde{v}_2 + \eta g - \sigma\Delta\eta = \tilde{B}(\tilde{\boldsymbol{v}}, \eta) \qquad \text{on } z = 0. \tag{2.4e}$$

Here

$$\tilde{\boldsymbol{L}}(\eta) = \begin{pmatrix} \frac{\eta}{d}U_1 \\ \frac{\eta}{d}U_2 \\ -\eta_x\frac{z+d}{d}U_1 - \eta_y\frac{z+d}{d}U_2 \end{pmatrix}, \quad \tilde{\boldsymbol{N}}(\tilde{\boldsymbol{v}}, \eta) = \boldsymbol{N}(\boldsymbol{U} + \tilde{\boldsymbol{v}}, \eta) - D\boldsymbol{N}[\boldsymbol{U}, 0](\tilde{\boldsymbol{v}}, \eta)$$

and

$$\tilde{B}(\tilde{\boldsymbol{v}}, \eta) = -\frac{1}{2}[B(\boldsymbol{U} + \tilde{\boldsymbol{v}}, \eta) - DB[\boldsymbol{U}, 0](\tilde{\boldsymbol{v}}, \eta) - Q(c_1, c_2)]_{z=0} + \sigma\nabla \cdot \left( \frac{\nabla\eta|\nabla\eta|^2}{\sqrt{1+|\nabla\eta|^2}(\sqrt{1+|\nabla\eta|^2}+1)} \right).$$

Note that $\tilde{\boldsymbol{N}}(\tilde{\boldsymbol{v}}, 0) = 0$ and $D\tilde{\boldsymbol{N}}[0, 0] = 0$. Similarly, $\tilde{B}(0, 0) = 0$ (since $B(\boldsymbol{U}, 0) = Q(c_1, c_2)$) and $D\tilde{B}[0, 0] = 0$.

We can simplify the linear part of problem (2.4) by introducing the new variable

$$\boldsymbol{v} = \tilde{\boldsymbol{v}} - \boldsymbol{w}^\eta - \tilde{\boldsymbol{U}}^\eta, \text{ where } \boldsymbol{w}^\eta = \begin{pmatrix} \frac{\eta}{d}U_1 + \alpha\eta\frac{z+d}{d}U_2 \\ \frac{\eta}{d}U_2 - \alpha\eta\frac{z+d}{d}U_1 \\ -\eta_x\frac{z+d}{d}U_1 - \eta_y\frac{z+d}{d}U_2 \end{pmatrix},$$

and $\tilde{\boldsymbol{U}}^\eta = \boldsymbol{U}[\tilde{c}_1^\eta, \tilde{c}_2^\eta]$ with

$$\begin{aligned}\int_{\Omega_{00}^0} \tilde{U}_j^\eta \,\mathrm{d}V &= \int_{\Omega_{00}^\eta} U_j \,\mathrm{d}V - \int_{\Omega_{00}^0} U_j \,\mathrm{d}V - \int_{\Omega_{00}^0} w_j^\eta \,\mathrm{d}V \\ &= \int_{\Omega_{00}^0} U_j \left(z + \eta(\boldsymbol{x}')\left(\frac{z}{d}+1\right)\right) \frac{\eta(\boldsymbol{x}')+d}{d} \,\mathrm{d}V \\ &\quad - \int_{\Omega_{00}^0} U_j(z) \frac{\eta(\boldsymbol{x}')+d}{d} \,\mathrm{d}V \\ &\quad - \int_{\Omega_{00}^0} U_j'(z) \frac{\eta(\boldsymbol{x}')(z+d)}{d} \,\mathrm{d}V, \quad j = 1, 2.\end{aligned}$$

This is a linear system of equations for $\tilde{c}_1^\eta$, $\tilde{c}_2^\eta$, which is uniquely solvable if and only if $\alpha d \notin 2\pi\mathbb{Z} \setminus \{0\}$, since

$$\int_{\Omega_{00}^0} \tilde{U}_1^\eta \,\mathrm{d}V = \frac{\sin(\alpha d)}{\alpha} \tilde{c}_1^\eta + \frac{(\cos(\alpha d) - 1)}{\alpha} \tilde{c}_2^\eta$$

and

$$\int_{\Omega_{00}^0} \tilde{U}_2^\eta \,\mathrm{d}V = -\frac{(\cos(\alpha d) - 1)}{\alpha} \tilde{c}_1^\eta + \frac{\sin(\alpha d)}{\alpha} \tilde{c}_2^\eta,$$

with obvious modifications if $\alpha = 0$. We shall make this assumption from now on. The transformation then gives

$$\begin{aligned} \nabla \times \boldsymbol{v} - \alpha \boldsymbol{v} &= \boldsymbol{G}(\boldsymbol{v}, \eta) && \text{in } \Omega^0, && \text{(2.5a)} \\ \nabla \cdot \boldsymbol{v} &= 0 && \text{in } \Omega^0, && \text{(2.5b)} \\ v_3 &= c_1 \eta_x + c_2 \eta_y && \text{on } z = 0, && \text{(2.5c)} \\ v_3 &= 0 && \text{on } z = -d, && \text{(2.5d)} \\ \int_{\Omega_{00}^0} v_j \,\mathrm{d}V &= 0 && j = 1, 2, && \text{(2.5e)} \\ c_1 v_1 + c_2 v_2 + g\eta - \sigma \Delta \eta &= R(\boldsymbol{v}, \eta) && \text{on } z = 0. && \text{(2.5f)} \end{aligned}$$

Here $\boldsymbol{G}(\boldsymbol{v}, \eta) = \nabla \times \tilde{\boldsymbol{N}}(\boldsymbol{v} + \boldsymbol{w}^\eta + \tilde{\boldsymbol{U}}^\eta, \eta)$ and $R(\boldsymbol{v}, \eta) = \tilde{B}(\boldsymbol{v} + \boldsymbol{w}^\eta + \tilde{\boldsymbol{U}}^\eta, \eta) - c_1 \tilde{c}_1^\eta - c_2 \tilde{c}_2^\eta$. Note that $\boldsymbol{G}$ is affine linear in its first argument. We have $\boldsymbol{G}(\boldsymbol{v}, 0) = 0$ and $D\boldsymbol{G}[0, 0] = 0$ as well as $R(0, 0) = 0$ and $DR[0, 0] = 0$. Therefore, the linearisation of (2.5) is the same as the formal linearisation of (1.6).

### *2.3. Reduction to the Surface*

We now go on to show that problem (1.6) (or, equivalently, problem (2.5)) can be reduced to a nonlinear pseudodifferential equation for the surface elevation $\eta$ in a neighbourhood of a laminar flow. To do this, we eliminate the vector field $\boldsymbol{v}$ from equation (2.5f) by solving (2.5a)–(2.5e) for $\boldsymbol{v}$. The solution $\boldsymbol{v} = \boldsymbol{v}(\eta, \boldsymbol{c})$ is expressed as an analytic operator of $\eta$ and $\boldsymbol{c}$ with $\boldsymbol{v}(0, \boldsymbol{c}) = 0$. Substituting this into (2.5f), we can rewrite problem (2.5) as the single equation

$$c_1[v_1(\eta, \boldsymbol{c})]_{z=0} + c_2[v_2(\eta, \boldsymbol{c})]_{z=0} + g\eta - \sigma\Delta\eta - R(\boldsymbol{v}(\eta, \boldsymbol{c}), \eta) = 0. \tag{2.6}$$

In order for the procedure to work, we need to impose the non-resonance condition

$$\sqrt{\alpha^2 - |\boldsymbol{k}|^2} \notin \frac{\pi}{d}\mathbb{Z}_+ \quad \text{for all } \boldsymbol{k} \in \Lambda' \text{ such that } |\boldsymbol{k}| < |\alpha|, \tag{2.7}$$

where $\mathbb{Z}_+ = \{1, 2, 3, \ldots\}$.

**Theorem 2.1.** *Let $d > 0$, $\alpha \in \mathbb{R}$ and $\Lambda'$ be given and assume that condition (2.7) holds. There exists a constant $r_0 > 0$ such that for any $\eta \in C^{2,\gamma}_{per,e}(\mathbb{R}^2)$ and any $\boldsymbol{c} \in \mathbb{R}^2$ problem (2.5a)–(2.5e) admits a unique solution $\boldsymbol{v} \in (C^{1,\gamma}_{per,e}(\overline{\Omega^0}))^2 \times C^{1,\gamma}_{per,o}(\overline{\Omega^0})$ provided $\|\eta\|_{C^{2,\gamma}(\mathbb{R}^2)} < r_0$. The constant $r_0$ depends only on $\alpha$ and $d$. Furthermore, the vector field $\boldsymbol{v}$ depends analytically on $\eta$ and $\boldsymbol{c}$. If $\eta$ is constant in the direction $\boldsymbol{\lambda}_j$, then so is $\boldsymbol{v}$.*

The proof is based on a perturbative approach. We first define suitable Banach spaces and prove a technical lemma for the unperturbed problem. After that we deal with the perturbed problem, thus proving the theorem.

In the analysis of (2.5a)–(2.5e) we assume that $\boldsymbol{v} \in Y$, where

$$Y := \left\{ \boldsymbol{v} \in (C^{1,\gamma}_{per,e}(\overline{\Omega^0}))^2 \times C^{1,\gamma}_{per,o}(\overline{\Omega^0}) : \nabla \cdot \boldsymbol{v} = 0 \text{ in } \Omega^0,\ v_3|_{\partial\Omega^{0,b}} = 0,\ \int_{\Omega^0_{00}} v_j \, \mathrm{d}V = 0,\ j = 1, 2 \right\}.$$

The corresponding range space for the operator

$$C_\alpha \colon \boldsymbol{v} \mapsto \left(\nabla \times \boldsymbol{v} - \alpha\boldsymbol{v}, v_3|_{\partial\Omega^{0,s}}\right)$$

is given by

$$Z := \left\{ (\boldsymbol{w}, f) \in ((C^{0,\gamma}_{per,e}(\overline{\Omega^0}))^2 \times C^{0,\gamma}_{per,o}(\overline{\Omega^0})) \times C^{1,\gamma}_{per,o}(\mathbb{R}^2) : \nabla \cdot \boldsymbol{w} = 0 \text{ in } \Omega^0 \right\}.$$

**Lemma 2.2.**

*(i) The operator $C_0 \colon Y \to Z$ is a linear isomorphism.*
*(ii) The operator $C_\alpha \colon Y \to Z$ is Fredholm of index zero for all $\alpha \in \mathbb{R}$.*
*(iii) If $\sqrt{\alpha^2 - |\boldsymbol{k}|^2} \notin \frac{\pi}{d}\mathbb{Z}_+$ for all $\boldsymbol{k} \in \Lambda'$, then $C_\alpha \colon Y \to Z$ is a linear isomorphism.*

**Proof.** Let us prove the first claim. We start with the injectivity. For a given $\boldsymbol{v} \in Y$ such that

$$\nabla \times \boldsymbol{v} = 0 \text{ in } \Omega^0, \tag{2.8}$$

and

$$v_3 = 0 \text{ on } \partial\Omega^{0,s},$$

we find that

$$\Delta v_3 = 0 \text{ in } \Omega^0,$$

while $v_3 = 0$ on $\partial\Omega^0$. Thus, $v_3 = 0$ identically. Using this fact and equation (2.8) restricted to the boundary, we obtain

$$\partial_z v_1 = \partial_z v_2 = 0 \text{ on } \partial\Omega^0.$$

Thus, the components $v_1$ and $v_2$ must be constant throughout $\Omega^0$ as they solve $\Delta v_1 = \Delta v_2 = 0$ with homogeneous Neumann conditions. But the integral assumptions in the definition of the space $Y$ require these constants to be zero and so $\boldsymbol{v}$ vanishes everywhere in $\Omega^0$ and we obtain the injectivity.

Next we turn to the surjectivity. We need to solve the equations

$$\nabla \times \boldsymbol{v} = \boldsymbol{w} \text{ in } \Omega^0, \tag{2.9a}$$

and

$$v_3 = f \text{ on } \partial\Omega^{0,s} \tag{2.9b}$$

for $(\boldsymbol{w}, f) \in Z$ with $\boldsymbol{v} \in Y$. Classical elliptic theory (see for instance Agmon, Douglis & Nirenberg [2, Theorem 6.30]) provides the existence of solutions $A_j \in C^{2,\gamma}_{per}(\overline{\Omega^0})$ to the equations

$$\begin{aligned} \Delta A_j &= w_j \text{ in } \Omega^0, \quad j = 1, 2, 3, \\ A_1 = A_2 &= \partial_z A_3 = 0 \text{ on } \partial\Omega^0 \end{aligned}$$

(note that integral of $w_3$ over $\Omega^0_{00}$ vanishes due to oddness). Furthermore, it's easily seen that $A_1$ and $A_2$ are even in $\boldsymbol{x}'$ and that $A_3$ can be chosen odd. We put

$$A = \nabla \cdot \boldsymbol{A}, \quad \boldsymbol{A} = (A_1, A_2, A_3)$$

and note that $A \in C^{1,\gamma}_{per}(\overline{\Omega^0})$ solves

$$\Delta A = 0 \text{ in } \Omega^0, \quad A = 0 \text{ on } \partial\Omega^0.$$

The unique solvability of the Dirichlet problem implies $\nabla \cdot \boldsymbol{A} = 0$ everywhere in $\Omega^0$. Moreover, we let $\varphi \in C^{2,\gamma}_{per}(\overline{\Omega^0})$ be the unique odd solution to

$$\begin{aligned} \Delta\varphi &= 0 \text{ in } \Omega^0, \\ \partial_z\varphi &= f \text{ on } \partial\Omega^{0,s}, \\ \partial_z\varphi &= 0 \text{ on } \partial\Omega^{0,b}, \end{aligned}$$

and set

$$\boldsymbol{v} = -\nabla \times \boldsymbol{A} + \nabla\varphi.$$

It is straightforward to verify that $\boldsymbol{v}$ satisfies (2.9a). Furthermore, the boundary conditions (2.9b) and

$$v_3 = 0 \text{ on } \partial\Omega^{0,b}$$

follow from the relation

$$v_3 = -\partial_x A_2 + \partial_y A_1 + \partial_z \varphi.$$

On the other hand, we find that

$$v_1 = -\partial_y A_3 + \partial_z A_2 + \partial_x \varphi, \quad v_2 = -\partial_z A_1 + \partial_x A_3 + \partial_y \varphi.$$

Now because $A_1 = A_2 = 0$ on the boundary, we find that

$$\int_{\Omega_{00}^0} v_j \,\mathrm{d}V = 0, \quad j = 1, 2.$$

Finally, the formulas also reveal that $v_1$ and $v_2$ are even in $\boldsymbol{x}'$, while $v_3$ is odd. Thus, $\boldsymbol{v} \in Y$ and the surjectivity is verified. To complete our proof of the first statement we use the inverse mapping theorem which ensures that $C_0 \colon Y \to Z$ is an isomorphism.

The second claim follows from an observation that the composition

$$C_0^{-1} \circ C_\alpha \colon Y \to Y$$

is a sum of the identity and a compact operator. Thus it is Fredholm of index zero and so is $C_\alpha$.

The last statement (iii) follows from (ii) since under the given assumption the kernel of $C_\alpha$ is trivial. Indeed, if $C_\alpha(\boldsymbol{v}) = 0$ for some $\boldsymbol{v} \in X$, then the Fourier coefficients solve

$$(\hat{v}_j^{(\boldsymbol{k})})'' + (\alpha^2 - |\boldsymbol{k}|^2)\hat{v}_j^{(\boldsymbol{k})} = 0 \ \text{ in } \ \Omega^0$$

for all $\boldsymbol{k} \in \Lambda'$ and $j = 1, 2, 3$. The component $\hat{v}_3^{(\boldsymbol{k})}$ also satisfies homogeneous Dirichlet boundary conditions both at $z = -d$ and $z = 0$. According to the assumption of the claim, we see that $\alpha^2 - |\boldsymbol{k}|^2$ is not an eigenvalue and hence $\hat{v}_3^{(\boldsymbol{k})}$ must be zero everywhere in $\Omega^0$. Using this fact, we can compute the third coordinate of $\nabla \times \boldsymbol{v} - \alpha \boldsymbol{v}$ to find

$$k_1 \hat{v}_2^{(\boldsymbol{k})} - k_2 \hat{v}_1^{(\boldsymbol{k})} = 0$$

for all $z \in [-d, 0]$, where $k_j = \boldsymbol{k} \cdot \boldsymbol{e}_j$, $j = 1, 2$. Moreover, since $\boldsymbol{v}$ is divergence free, we have

$$k_1 \hat{v}_1^{(\boldsymbol{k})} + k_2 \hat{v}_2^{(\boldsymbol{k})} = 0.$$

Taken together, this shows that

$$|\boldsymbol{k}|^2 \hat{v}_j^{(\boldsymbol{k})} = 0, \quad j = 1, 2.$$

Thus, $\boldsymbol{v} = 0$ identically and the kernel is trivial. This finishes the proof of the lemma. □

We are now ready to treat the perturbed problem.

**Proof of Theorem 2.1.** Note that $C_\alpha\colon Y \to Z$ is an isomorphism by the hypotheses of the theorem and Lemma 2.2. We can therefore rewrite problem (2.5a)–(2.5e) as

$$\boldsymbol{v} - C_\alpha^{-1}(\boldsymbol{G}(\boldsymbol{v}, \eta) - \boldsymbol{G}(0, \eta), 0) = C_\alpha^{-1}(\boldsymbol{G}(0, \eta), c_1\eta_x + c_2\eta_y),$$

where the left-hand side is a bounded linear operator on $Y$ (since $\boldsymbol{v} \mapsto \boldsymbol{G}(\boldsymbol{v}, \eta)$ is affine linear), which is analytic in $\eta \in U := \{\eta \in C^{2,\gamma}_{per}(\mathbb{R}^2) : \min \eta > -d\}$, and the right hand side is an analytic mapping $U \times \mathbb{R}^2 \to Y$. Since $\|C_\alpha^{-1}(\boldsymbol{G}(\boldsymbol{v}, \eta) - \boldsymbol{G}(0, \eta), 0)\|_Y = \mathcal{O}(\|\eta\|_{C^{2,\gamma}(\mathbb{R}^2)})\|\boldsymbol{v}\|_Y$ it follows by the analytic implicit-function theorem that the equation has a unique solution $\boldsymbol{v} = \boldsymbol{v}(\eta, \boldsymbol{c})$ which is analytic in $\eta$ and $\boldsymbol{c}$ for $\|\eta\|_{C^{2,\gamma}(\mathbb{R}^2)} < r_0$ with $r_0$ sufficiently small. To prove the last claim, we repeat the analysis in the subspace of $Y$ consisting of vector fields which are constant in the direction $\boldsymbol{\lambda}_j$ , noting that $C_\alpha$ and $\boldsymbol{G}$ preserve this property if $\eta$ also shares it. □

## 3. Analysis of the Linearised Problem

### *3.1. Dispersion Equation*

In this section we analyse the linearisation of (2.5). We choose to work with this problem rather than the reduced equation (2.6) since it is slightly more general (see condition (2.7)) and since it also gives direct information about the velocity field. However, it is of course straightforward to draw conclusions concerning the linearisation of problem (2.6) from this analysis. We aim to show that for a broad range of parameters the kernel of the linearised problem is exactly four-dimensional (and therefore two-dimensional when restricting to solutions satisfying the symmetry conditions (1.5)). For the convenience of the reader, the results are summarised in Proposition 3.1 at the end of the section.

The kernel of the linearisation of problem (2.5) is described by the system

$$\nabla \times \boldsymbol{v} - \alpha \boldsymbol{v} = 0 \qquad \text{in } \Omega^0, \tag{3.1a}$$

$$\nabla \cdot \boldsymbol{v} = 0 \qquad \text{in } \Omega^0, \tag{3.1b}$$

$$v_3 = c_1\eta_x + c_2\eta_y \qquad \text{on } z = 0, \tag{3.1c}$$

$$v_3 = 0 \qquad \text{on } z = -d, \tag{3.1d}$$

$$\int_{\Omega^0_{00}} v_j \,\mathrm{d}V = 0, \qquad \text{for } j = 1, 2, \tag{3.1e}$$

$$c_1 v_1 + c_2 v_2 + g\eta - \sigma\Delta\eta = 0 \qquad \text{on } z = 0. \tag{3.1f}$$

By Fourier analysis, it is enough to consider the Ansatz

$$\eta = \hat{\eta} e^{i\boldsymbol{k}\cdot\boldsymbol{x}'}, \quad \boldsymbol{v} = (\hat{v}_1(z), \hat{v}_2(z), \hat{v}_3(z))e^{i\boldsymbol{k}\cdot\boldsymbol{x}'}$$

with $\boldsymbol{k} = (k, l) \in \Lambda'$ in order to find periodic solutions of these equations. We split the analysis into four cases, in which $\sqrt{\cdot}$ denotes the principal branch of the square root.

**Case I:** $\sqrt{\alpha^2 - |\boldsymbol{k}|^2} \notin \frac{\pi}{d}\mathbb{Z}_+$ and $\boldsymbol{k} \neq (0,0)$.

The first two equations (3.1a) and (3.1b) imply that

$$\hat{v}_3''(z) + (\alpha^2 - |\boldsymbol{k}|^2)\hat{v}_3(z) = 0, \quad z \in (-d, 0). \tag{3.2}$$

Taking into account the boundary conditions (3.1c) and (3.1d), we find

$$\hat{v}_3(z) = \lambda\hat{\eta}\phi(z),$$

where $\lambda = i(\boldsymbol{c} \cdot \boldsymbol{k})$ and

$$\phi(z) = \begin{cases} \frac{\sin(\sqrt{\alpha^2 - |\boldsymbol{k}|^2}(z+d))}{\sin(\sqrt{\alpha^2 - |\boldsymbol{k}|^2}d)} & \text{if } |\alpha| > |\boldsymbol{k}|, \\ \frac{z+d}{d} & \text{if } |\alpha| = |\boldsymbol{k}|, \\ \frac{\sinh(\sqrt{|\boldsymbol{k}|^2 - \alpha^2}(z+d))}{\sinh(\sqrt{|\boldsymbol{k}|^2 - \alpha^2}d)} & \text{if } |\alpha| < |\boldsymbol{k}|. \end{cases}$$

Now since $\boldsymbol{k} \neq (0,0)$, we obtain

$$\hat{v}_1(z) = i\lambda|\boldsymbol{k}|^{-2}\hat{\eta}(k\phi'(z) + \alpha l\phi(z)),$$
$$\hat{v}_2(z) = i\lambda|\boldsymbol{k}|^{-2}\hat{\eta}(l\phi'(z) - \alpha k\phi(z)).$$

Note that (3.1e) is satisfied automatically since $\boldsymbol{k} \neq (0,0)$. Substituting $\boldsymbol{v}$ into (3.1f) and assuming $\hat{\eta} \neq 0$ (otherwise $\boldsymbol{v}$ vanishes), we arrive at the dispersion equation

$$\rho(\boldsymbol{c}, \boldsymbol{k}) := g + \sigma|\boldsymbol{k}|^2 - \frac{(\boldsymbol{c} \cdot \boldsymbol{k})^2}{|\boldsymbol{k}|^2}\kappa(|\boldsymbol{k}|) + \alpha\frac{(\boldsymbol{c} \cdot \boldsymbol{k})(\boldsymbol{c} \cdot \boldsymbol{k}_\perp)}{|\boldsymbol{k}|^2} = 0, \tag{3.3}$$

where

$$\kappa(|\boldsymbol{k}|) := \phi'(0; |\boldsymbol{k}|) = \begin{cases} \sqrt{\alpha^2 - |\boldsymbol{k}|^2}\cot(\sqrt{\alpha^2 - |\boldsymbol{k}|^2}d) & \text{if } |\alpha| > |\boldsymbol{k}|, \\ \frac{1}{d} & \text{if } |\alpha| = |\boldsymbol{k}|, \\ \sqrt{|\boldsymbol{k}|^2 - \alpha^2}\coth(\sqrt{|\boldsymbol{k}|^2 - \alpha^2}d) & \text{if } |\alpha| < |\boldsymbol{k}| \end{cases} \tag{3.4}$$

and

$$\boldsymbol{k}_\perp = (-l, k).$$

This is an equation for $\boldsymbol{k}$ which will be analysed below.

**Case II:** $|\alpha| \notin \frac{\pi}{d}\mathbb{Z}_+$ and $\boldsymbol{k} = (0,0)$.

This corresponds to the constant function $\eta = \hat{\eta}$. The vector field $\boldsymbol{v}$ solving (3.1a)–(3.1d) with no $\boldsymbol{x}'$-dependence must coincide with a laminar flow $\boldsymbol{U}[\tilde{c}_1, \tilde{c}_2]$ for some constants $\tilde{c}_1, \tilde{c}_2 \in \mathbb{R}$, but the condition (3.1e) forces them to be zero. Finally, condition (3.1f) forces $\hat{\eta} = 0$. Thus, in this case we find no non-trivial solutions to the linearised problem.

**Case III:** $\sqrt{\alpha^2 - |\boldsymbol{k}|^2} \in \frac{\pi}{d}\mathbb{Z}_+$ and $\boldsymbol{k} \neq (0,0)$.

Just as in Case I, we find that $\hat{v}_3$ solves (3.2). The condition $\hat{v}_3(-d) = 0$ is implied by (3.1d) and, since $\sqrt{\alpha^2 - |\boldsymbol{k}|^2} \in \frac{\pi}{d}\mathbb{Z}_+$, we necessarily have $\hat{v}_3(0) = 0$. Thus, $\hat{v}_3(z) = \lambda\phi_0(z)$, where $\lambda$ is an arbitrary constant and $\phi_0(z) = \sin(\pi nz/d)$, $n \in \mathbb{Z}_+$, solves (3.2) with homogeneous Dirichlet boundary conditions. Furthermore, (3.1c) reduces to $(\boldsymbol{c}\cdot\boldsymbol{k})\hat{\eta} = 0$. Thus, either $\hat{\eta}$ or $\boldsymbol{c}\cdot\boldsymbol{k}$ vanishes. The functions $\hat{v}_1$ and $\hat{v}_2$ are given by the same formulas as in Case I but with the function $\phi$ replaced by $\phi_0$ and without $\hat{\eta}$. Again (3.1e) follows from $\boldsymbol{k} \neq (0, 0)$. Finally, relation (3.1f) leads to

$$\frac{i\lambda(\boldsymbol{c}\cdot\boldsymbol{k})}{|\boldsymbol{k}|^2}\phi_0'(0) + \hat{\eta}[g + \sigma|\boldsymbol{k}|^2] = 0.$$

If $\boldsymbol{c}\cdot\boldsymbol{k} = 0$, we see from this that $\hat{\eta} = 0$ while there are no restrictions on $\lambda$. If on the other hand, $\boldsymbol{c}\cdot\boldsymbol{k} \neq 0$, we saw before that $\hat{\eta} = 0$ and this forces $\lambda = 0$ (note that $\phi_0'(0) \neq 0$) and hence $\boldsymbol{v} = 0$. Thus $\eta = 0$, but $\boldsymbol{v}$ need not vanish if $\boldsymbol{c}\cdot\boldsymbol{k} = 0$.

**Case IV:** $|\alpha| \in \frac{\pi}{d}\mathbb{Z}_+$ and $\boldsymbol{k} = (0, 0)$.

As in Case II, we get that $\eta = \hat{\eta}$ is constant and $\boldsymbol{v} = \boldsymbol{U}[\tilde{c}_1, \tilde{c}_2]$ for some constants $\tilde{c}_1, \tilde{c}_2 \in \mathbb{R}$. If $|\alpha|$ is an odd multiple of $\pi/d$, (3.1e) forces $\tilde{c}_1 = \tilde{c}_2 = 0$ and then (3.1f) leads to $\hat{\eta} = 0$. However, in the even case $\tilde{c}_1$ and $\tilde{c}_2$ are arbitrary and (3.1f) leads to

$$\hat{\eta} = -\frac{\boldsymbol{c}\cdot\tilde{\boldsymbol{c}}}{g}.$$

The last two cases are included for completeness and future reference, but we will avoid them in the further analysis by assuming the non-resonance condition (2.7). A complete analysis of equation (3.3) is a complicated problem. Our aim here is to find sufficient conditions on the problem parameters that guarantee that there exists some $\boldsymbol{c} = \boldsymbol{c}^\star$ such that the dispersion equation $\rho(\boldsymbol{c}^\star, \boldsymbol{k}) = 0$ has exactly four different roots in the dual lattice $\Lambda'$. Note that the roots come in pairs $\pm\boldsymbol{k}$, so that the dimension of the solution space is halved when we consider solutions with the symmetries (1.5). We use a geometric approach and restrict ourselves to the case when the roots are generators of $\Lambda'$. We will also restrict attention to the case $\alpha \neq 0$ in the main part of the analysis and leave the irrotational case to Remarks 3.2, 3.5 and 4.2 (see also the references mentioned in the introduction). Let us assume that the constants $\alpha$, $\sigma$ and $d$ are fixed. Then for a given $\boldsymbol{k} \neq \boldsymbol{0}$ we want to describe the set of all $\boldsymbol{c} \in \mathbb{R}^2$ such that (3.3) holds true. In other words, we are going to analyse the zero level set

$$\rho(\boldsymbol{c}, \boldsymbol{k}) = 0, \tag{3.5}$$

where the vector $\boldsymbol{k}$ is fixed. For this purpose we put

$$x = \frac{\boldsymbol{c}\cdot\boldsymbol{k}}{|\boldsymbol{k}|}, \quad y = \frac{\boldsymbol{c}\cdot\boldsymbol{k}_\perp}{|\boldsymbol{k}|}$$

and write equation (3.5) in the form

$$\kappa(|\boldsymbol{k}|)x^2 = a(|\boldsymbol{k}|) + \alpha xy, \qquad a(|\boldsymbol{k}|) := g + \sigma|\boldsymbol{k}|^2.$$

We can solve for $y$ and get a curve of solutions in the $xy$-plane given by

$$y = \frac{\kappa(|\boldsymbol{k}|)}{\alpha}x - \frac{a(|\boldsymbol{k}|)}{\alpha x}.$$

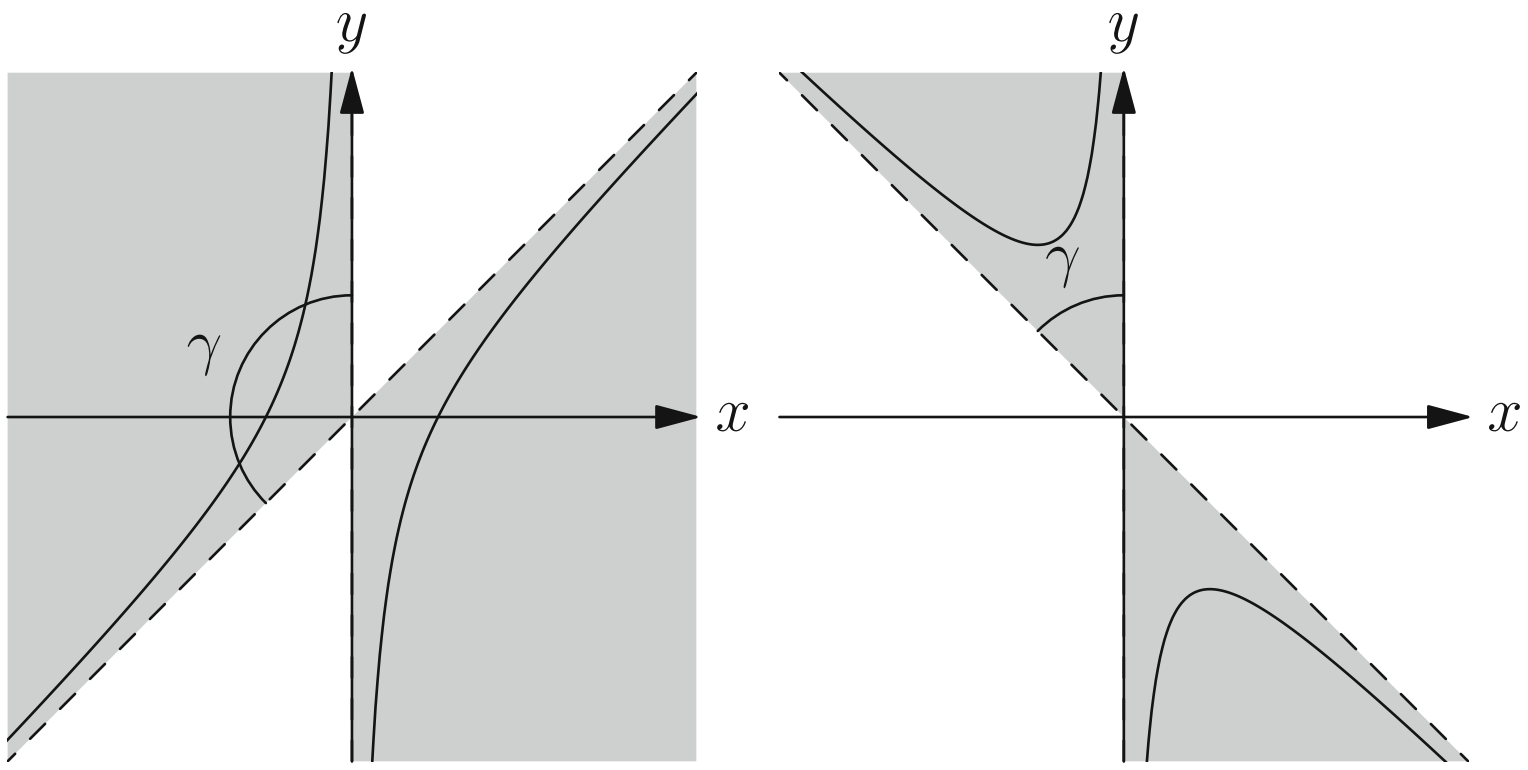

**(a)** $\kappa(|\boldsymbol{k}|) > 0$ **(b)** $\kappa(|\boldsymbol{k}|) < 0$

**Fig. 4.** The hyperbolas are solutions to the dispersion equation (3.3) in the $xy$-plane for $\alpha > 0$

The curve can be recognised as a hyperbola with the $y$-axis as one asymptote. We denote by

$$\gamma = \frac{\pi}{2} + \arctan\left(\frac{\kappa(|\boldsymbol{k}|)}{|\alpha|}\right)$$

the *angle between the asymptotes* of one branch of the hyperbola (Figure 4). It is clear that $\gamma \in (0, \pi)$ and that $\gamma$ is obtuse if $\kappa(|\boldsymbol{k}|)$ is positive and acute if it is negative. We call the open set $\{(x, y) : \kappa(|\boldsymbol{k}|)x^2 > \alpha xy\}$, which contains the hyperbola and is delimited by its asymptotes, the *set between the asymptotes*. It is the shaded set in Figures 4a and 4b.

To express this curve in $(c_1, c_2)$ coordinates we note that $x = c_1 \cos(\theta) + c_2 \sin(\theta)$ and $y = -c_1 \sin(\theta) + c_2 \cos(\theta)$, where $\theta$ is the angle that $\boldsymbol{k}$ makes with $\boldsymbol{e}_1$. Hence

$$\begin{pmatrix} c_1 \\ c_2 \end{pmatrix} = \begin{pmatrix} \cos(\theta) & -\sin(\theta) \\ \sin(\theta) & \cos(\theta) \end{pmatrix} \begin{pmatrix} x \\ y \end{pmatrix},$$

so going from $(x, y)$ to $(c_1, c_2)$ is a counterclockwise rotation by the angle $\theta$. We denote this curve of solutions in the $(c_1, c_2)$ plane by $\mathcal{C}(\boldsymbol{k})$. Note that the equation $\rho(\boldsymbol{c}, \boldsymbol{k}_j) = 0$ can also be written in the form

$$\boldsymbol{c}^T A_j \boldsymbol{c} = 1, \quad j = 1, 2, \tag{3.6}$$

where $A_j$ are real indefinite symmetric $2 \times 2$ matrices. The set between the asymptotes of $\mathcal{C}(\boldsymbol{k}_j)$ is therefore given by $\{\boldsymbol{c} : \boldsymbol{c}^T A_j \boldsymbol{c} > 0\}$. Now we want to find linearly independent vectors $\boldsymbol{k}_1$ and $\boldsymbol{k}_2$ so that there is a point of intersection of $\mathcal{C}(\boldsymbol{k}_1)$ and $\mathcal{C}(\boldsymbol{k}_2)$. Clearly a sufficient condition is that the sets between the asymptotes of $\mathcal{C}(\boldsymbol{k}_1)$ and $\mathcal{C}(\boldsymbol{k}_2)$ have nonempty intersection, but one is not contained in the other. (3.7)

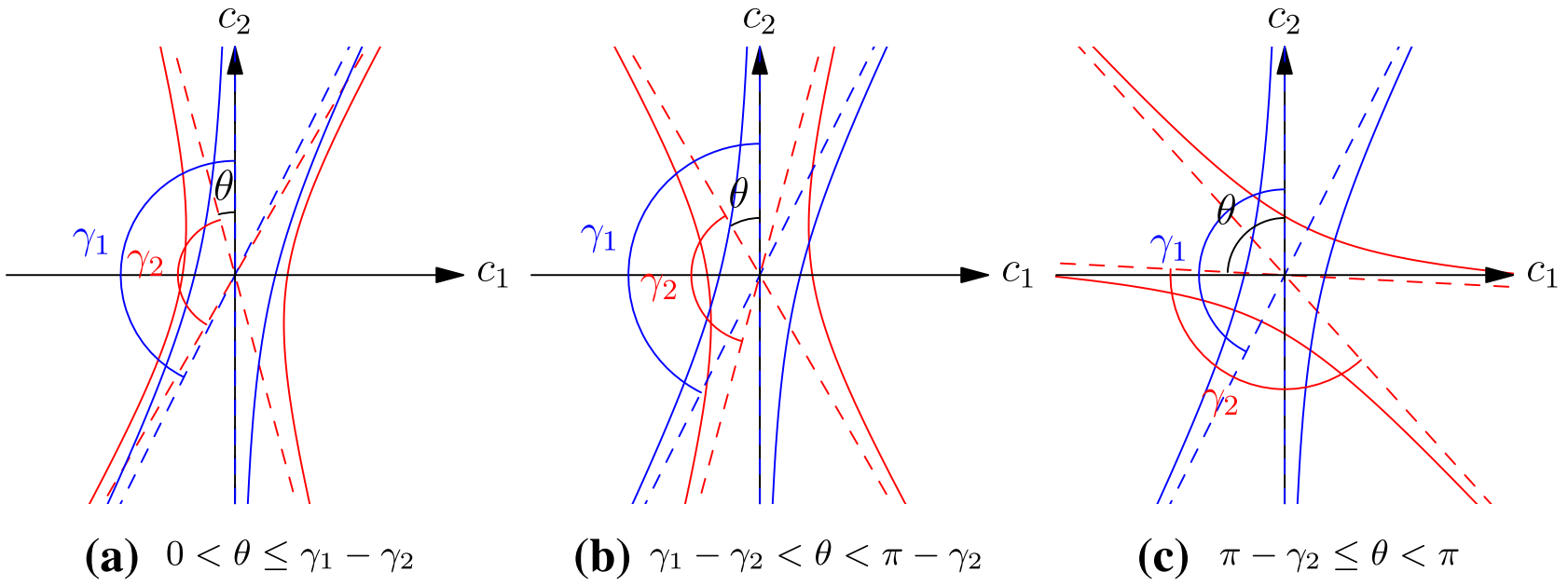

**Fig. 5.** Intersection points of the hyperbolas correspond to common solutions $\boldsymbol{c} = (c_1, c_2)$ of the dispersion equation (3.3) for two different wave vectors $\boldsymbol{k}_1$ and $\boldsymbol{k}_2$. The figures illustrate the sufficient condition (3.8) for an intersection in the case $\alpha, \kappa(|\boldsymbol{k}_1|), \kappa(|\boldsymbol{k}_2|) > 0$

We now analyse this in more detail in the special case when $\alpha$, $\kappa(|\boldsymbol{k}_1|)$ and $\kappa(|\boldsymbol{k}_2|)$ are all positive. In that case we get two hyperbolas of the kind shown in Figure 4a. Note that if $\alpha \neq 0$ we we can always assume that it is positive by interchanging $x$ and $y$, and $u_1$ and $u_2$. Without loss of generality we can assume that $\boldsymbol{k}_1$ is parallel with $\boldsymbol{e}_1$ and that $\boldsymbol{k}_2$ makes an angle $\theta$ with $\boldsymbol{k}_1$. Moreover we can always choose the generators in such a way that $0 < \theta < \pi$ by changing $\boldsymbol{k}_2$ to $-\boldsymbol{k}_2$ if necessary. By possibly relabelling, we may assume that $\gamma_1 \geq \gamma_2$ where $\gamma_i$ is the angle between the asymptotes of $\mathcal{C}(\boldsymbol{k}_i)$. We see in Figure 5a that if $0 < \theta \leq \gamma_1 - \gamma_2$ then the set between the asymptotes of $\mathcal{C}(\boldsymbol{k}_2)$ is completely contained in the set between the asymptotes of $\mathcal{C}(\boldsymbol{k}_1)$; in Figures 5b and 5c we see that if $\gamma_1 - \gamma_2 < \theta < \pi$ then condition (3.7) is satisfied (note that $\gamma_1 - \gamma_2 < \pi - \gamma_2$). In summary we necessarily get intersection between $\mathcal{C}(\boldsymbol{k}_1)$ and $\mathcal{C}(\boldsymbol{k}_2)$ if

$$\gamma_1 - \gamma_2 < \theta < \pi. \tag{3.8}$$

In the subsequent analysis it will be important to make sure that the only solutions to the dispersion equation (3.5) for a fixed $\boldsymbol{c} = \boldsymbol{c}^\star$ are the generators $\pm\boldsymbol{k}_1$ and $\pm\boldsymbol{k}_2$. This property is expected to hold for generic values of the parameters since three hyperbolas in the plane generally don't have a common point of intersection. However, verifying it analytically is non-trivial. We content ourselves with analysing the case of a symmetric lattice, $|\boldsymbol{k}_1| = |\boldsymbol{k}_2| = k > 0$. In that case $\gamma_1 = \gamma_2$ and as before we assume that the angles are obtuse (meaning that $\kappa(k) = \kappa(|\boldsymbol{k}_1|) = \kappa(|\boldsymbol{k}_2|) > 0$) and that $\alpha > 0$. Condition (3.8) is clearly satisfied.

It's convenient to assume that the generators have the form

$$\boldsymbol{k}_1 = k(\cos\omega, \sin\omega), \quad \boldsymbol{k}_2 = k(\cos\omega, -\sin\omega),$$

with $\omega \in (0, \pi/2)$, which can always be achieved after rotating and relabelling. Note that the angle between $\boldsymbol{k}_1$ and $\boldsymbol{k}_2$ is $2\omega$. Similarly, we write

$$\boldsymbol{c} = c(\cos\varphi, \sin\varphi).$$

The dispersion equation $\rho(\boldsymbol{c}, \boldsymbol{k}_j) = 0$ then takes the form

$$c^2\Big(\kappa(k)(\cos\omega\cos\varphi + \sin\omega\sin\varphi)^2 \\ - \alpha(\cos\omega\cos\varphi + \sin\omega\sin\varphi)(\cos\omega\sin\varphi - \sin\omega\cos\varphi)\Big) = g + \sigma k^2$$

and

$$c^2\Big(\kappa(k)(\cos\omega\cos\varphi - \sin\omega\sin\varphi)^2 \\ - \alpha(\cos\omega\cos\varphi - \sin\omega\sin\varphi)(\cos\omega\sin\varphi + \sin\omega\cos\varphi)\Big) = g + \sigma k^2$$

for $j = 1$ and $j = 2$, respectively. It follows that $(c^\star)^2$ is proportional to $g + \sigma k^2$ (with proportionality constant independent of $g$ and $\sigma$) for any solution $\boldsymbol{c}^\star$.

Expanding the two equations we obtain

$$\begin{aligned} \nu &= \cos^2\omega(\kappa(k)\cos^2\varphi - \alpha\cos\varphi\sin\varphi) + \sin^2\omega(\kappa(k)\sin^2\varphi + \alpha\cos\varphi\sin\varphi) \\ &\quad + \cos\omega\sin\omega(2\kappa(k)\cos\varphi\sin\varphi + \alpha(\cos^2\varphi - \sin^2\varphi)) \\ &= \cos^2\omega(\kappa(k)\cos^2\varphi - \alpha\cos\varphi\sin\varphi) + \sin^2\omega(\kappa(k)\sin^2\varphi + \alpha\cos\varphi\sin\varphi) \\ &\quad - \cos\omega\sin\omega(2\kappa(k)\cos\varphi\sin\varphi + \alpha(\cos^2\varphi - \sin^2\varphi)), \end{aligned}$$

where

$$\nu = \frac{g + \sigma k^2}{(c^\star)^2}.$$

Taking the difference of the two different expressions for $\nu$ and using the fact that $\sin\omega\cos\omega \neq 0$, we obtain that

$$2\kappa(k)\cos\varphi\sin\varphi + \alpha(\cos^2\varphi - \sin^2\varphi) = 0$$

and hence

$$\tan(2\varphi) = -\frac{\alpha}{\kappa(k)}. \tag{3.9}$$

In particular, we note that $\boldsymbol{c}^\star$ is not parallel with either of the coordinate axes if $\alpha \neq 0$, in contrast to the irrotational case. We also note that we get

$$\begin{aligned} \nu &= \cos^2\omega(\kappa(k)\cos^2\varphi - \alpha\cos\varphi\sin\varphi) + \sin^2\omega(\kappa(k)\sin^2\varphi + \alpha\cos\varphi\sin\varphi) \\ &= (1 - 2\sin^2\omega)(\kappa(k)\cos^2\varphi - \alpha\cos\varphi\sin\varphi) + \kappa(k)\sin^2\omega. \end{aligned} \tag{3.10}$$

The dispersion equation $\rho(\boldsymbol{c}^\star, \boldsymbol{k}) = 0$ for a general lattice vector

$$\boldsymbol{k} = k((n_1 + n_2)\cos\omega, (n_1 - n_2)\sin\omega), \qquad n_1, n_2 \in \mathbb{Z},$$

can be written as

$$(c^\star)^2\Big(((n_1 + n_2)\cos\omega\cos\varphi + (n_1 - n_2)\sin\omega\sin\varphi)^2\kappa(kq_{n_1n_2}) \\ - \alpha((n_1 + n_2)\cos\omega\cos\varphi + (n_1 - n_2)\sin\omega\sin\varphi)((n_1 + n_2)\cos\omega\sin\varphi \\ - (n_1 - n_2)\sin\omega\cos\varphi)\Big) = (g + \sigma k^2 q_{n_1n_2}^2)q_{n_1n_2}^2,$$

where

$$q_{n_1 n_2} = \sqrt{(n_1 + n_2)^2 \cos^2 \omega + (n_1 - n_2)^2 \sin^2 \omega}.$$

Since the left hand side of the dispersion equation is proportional to $g + \sigma k^2$, we either get equality for all $\sigma$ or for at most one $\sigma > 0$. However, the former can happen only if $q_{n_1 n_2} = 1$ and the proportionality constant is 1. We have

$$q_{n_1 n_2}^2 \geq (|n_1| - |n_2|)^2$$

with equality only if either $n_1$ or $n_2$ vanishes or if both are non-zero and $\omega$ is an integer multiple of $\pi/2$. It follows that the only way to get $q_{n_1 n_2} = 1$ is if $(n_1, n_2) = (\pm 1, 0)$ or $(0, \pm 1)$, or if $n_1 = \pm n_2$. The former is the trivial case when $\boldsymbol{k} = \pm \boldsymbol{k}_j$, $j = 1, 2$. In the latter case, we get

$$q_{n_1 n_2}^2 = q_{n_1 n_1}^2 = 4 n_1^2 \cos^2 \omega = 1$$

or

$$q_{n_1 n_2}^2 = q_{n_1 (-n_1)}^2 = 4 n_1^2 \sin^2 \omega = 1.$$

We shall now show that this leads to a contradiction. For simplicity we restrict attention to the case $n_1 = n_2$, the other one being completely analogous. The dispersion equation simplifies to

$$\nu = \kappa(k) \cos^2 \varphi - \alpha \cos \varphi \sin \varphi. \tag{3.11}$$

Substituting this into (3.10), we get

$$\nu = (1 - 2 \sin^2 \omega) \nu + \kappa(k) \sin^2 \omega,$$

and hence,

$$\nu = \frac{\kappa(k)}{2}.$$

Substituting this result into (3.11), we get

$$\kappa(k)(2 \cos^2 \varphi - 1) = 2\alpha \cos \varphi \sin \varphi$$

or

$$\tan(2\varphi) = \frac{\kappa(k)}{\alpha},$$

which contradicts (3.9).

We can summarise the results of this section as follows:

**Proposition 3.1.**

*(i) Assume that the non-resonance condition* (2.7) *is satisfied for the dual lattice* $\Lambda'$. *Then the dimension of the space of solutions* $(\boldsymbol{v}, \eta) \in ((C^{1,\gamma}_{per,e}(\overline{\Omega^0}))^2 \times C^{1,\gamma}_{per,o}(\overline{\Omega^0})) \times C^{2,\gamma}_{per,e}(\mathbb{R}^2)$ *to the linearised problem* (3.1) *is equal to half the number of solutions* $\boldsymbol{k} \in \Lambda'$ *of the dispersion equation* (3.3) *(since solutions to the dispersion relation occur in pairs* $\pm\boldsymbol{k}$*).*

*(ii) If, in addition,* $\alpha \neq 0$ *and the sufficient condition* (3.7) *holds for the generators* $\boldsymbol{k}_1$ *and* $\boldsymbol{k}_2$ *of* $\Lambda'$, *then there exists a constant* $\boldsymbol{c}^\star$ *such that the solution space is at least two-dimensional when* $\boldsymbol{c} = \boldsymbol{c}^\star$.

*(iii) If* $\alpha$, $\kappa(|\boldsymbol{k}_1|)$ *and* $\kappa(|\boldsymbol{k}_2|)$ *are positive, where* $\kappa$ *is defined in* (3.4)*, then condition (3.7) is satisfied if the angle* $\theta$ *between the generators satisfies* (3.8).

*(iv) If, in addition to the conditions in (i) and (iii), the lattice is symmetric, that is,* $|\boldsymbol{k}_1| = |\boldsymbol{k}_2|$, *then the solution space is exactly two-dimensional for all but countably many values of* $\sigma$.

**Remark 3.2.** In the irrotational case, $\alpha = 0$, the hyperbolas instead become straight lines, $x = \pm\sqrt{a(|\boldsymbol{k}|)/\kappa(|\boldsymbol{k}|)}$. The rotated lines $\mathcal{C}(\boldsymbol{k}_1)$ and $\mathcal{C}(\boldsymbol{k}_2)$ will clearly intersect as soon as $\boldsymbol{k}_1$ and $\boldsymbol{k}_2$ are not parallel. The multiplicity analysis for symmetric lattices applies also in the irrotational case, but can then be simplified and expanded; see Section 7 of REEDER & SHINBROT [31].

### *3.2. Transversality Condition*

The local bifurcation theory that we are going to apply requires the bifurcation point (a laminar flow in our case) to satisfy a transversality condition. This can be stated as the condition that the vectors

$$\nabla_{\boldsymbol{c}}\rho(\boldsymbol{c}^\star, \boldsymbol{k}_1) \text{ and } \nabla_{\boldsymbol{c}}\rho(\boldsymbol{c}^\star, \boldsymbol{k}_2) \text{ are not parallel,} \tag{3.12}$$

where $\pm\boldsymbol{k}_1$ and $\pm\boldsymbol{k}_2$ are assumed to be only solutions to the equation $\rho(\boldsymbol{c}^\star, \boldsymbol{k}) = 0$. It turns out that (3.12) is automatically satisfied under the conditions discussed above.

**Proposition 3.3.** *Condition* (3.7) *is sufficient for the transversality condition* (3.12).

**Proof.** Assume instead that $\mathcal{C}(\boldsymbol{k}_1)$ and $\mathcal{C}(\boldsymbol{k}_2)$ intersect tangentially at $\boldsymbol{c}^\star$. Writing the equations $\rho(\boldsymbol{c}, \boldsymbol{k}_j) = 0$ in the form (3.6), we see that $A_1\boldsymbol{c}^\star = \lambda A_2\boldsymbol{c}^\star$ for some $\lambda \in \mathbb{R}$, and from $(\boldsymbol{c}^\star)^T A_j\boldsymbol{c}^\star = 1$, $j = 1, 2$, we get that $\lambda = 1$. Hence, $(A_1 - A_2)\boldsymbol{c}^\star = 0$. We now either get $A_1 \geq A_2$ or vice versa, depending on whether the other eigenvalue of the symmetric matrix $A_1 - A_2$ is nonnegative or nonpositive. For definiteness, we shall assume the former since it is consistent with the assumptions in the previous section. But this implies that the hyperbola $\boldsymbol{c}^T A_2\boldsymbol{c} = 1$ is contained in the set $\boldsymbol{c}^T A_1\boldsymbol{c} \geq 1$. Moreover, the set $\{\boldsymbol{c} : \boldsymbol{c}^T A_2\boldsymbol{c} > 0\}$ between the asymptotes of $\mathcal{C}(\boldsymbol{k}_2)$ is contained in the set $\{\boldsymbol{c} : \boldsymbol{c}^T A_1\boldsymbol{c} > 0\}$ between the asymptotes of $\mathcal{C}(\boldsymbol{k}_1)$, contradicting (3.7). □

**Remark 3.4.** In particular, condition (3.8) is sufficient for transversality in the special case $\alpha, \kappa(|\boldsymbol{k}_1|), \kappa(|\boldsymbol{k}_2|) > 0$ by the discussion in the previous section.

**Remark 3.5.** In the irrotational case, $\alpha = 0$, the transversality condition is automatically satisfied as soon as $\boldsymbol{k}_1$ and $\boldsymbol{k}_2$ are not parallel since the lines $\rho(\boldsymbol{c}^\star, \boldsymbol{k}_j) = 0$ then intersect transversally.

## 4. Main Result

Now we formulate our main theorem, providing the existence of three-dimensional steady gravity-capillary waves with vorticity.

**Theorem 4.1.** *Let $\alpha \in \mathbb{R}$, $\sigma > 0$ and the depth $d > 0$ be given, as well as a laminar flow $\boldsymbol{U}[c_1^\star, c_2^\star]$. Furthermore, let $\Lambda'$ be the dual lattice generated by the linearly independent vectors $\boldsymbol{k}_1, \boldsymbol{k}_2 \in \Lambda'$. Assume that*

*(i) the non-resonance condition (2.7) holds;*
*(ii) within the lattice $\Lambda'$, the dispersion equation (3.3) with $c_1 = c_1^\star, c_2 = c_2^\star$ has exactly four roots $\pm\boldsymbol{k}_1$ and $\pm\boldsymbol{k}_2$;*
*(iii) the transversality condition (3.12) holds.*

*Then there exists a neighbourhood of zero $W = B_\epsilon(0; \mathbb{R}^2) \subset \mathbb{R}^2$ and real-analytic functions $\delta_1, \delta_2 : W \to \mathbb{R}$ satisfying $\delta_1, \delta_2 = \mathcal{O}(|\boldsymbol{t}|^2)$ as $|\boldsymbol{t}| \to 0$, and such that for any $\boldsymbol{t} = (t_1, t_2) \in W$ there is a solution $(\boldsymbol{v}, \eta) \in ((C^{1,\gamma}_{per,e}(\overline{\Omega^0}))^2 \times C^{1,\gamma}_{per,o}(\overline{\Omega^0})) \times C^{2,\gamma}_{per,e}(\mathbb{R}^2)$ of problem (2.5) with*

$$c_1 = c_1^\star + \delta_1(\boldsymbol{t}), \quad c_2 = c_2^\star + \delta_2(\boldsymbol{t}), \quad Q = Q(c_1, c_2)$$

*such that*

$$\eta(\boldsymbol{x}') = t_1 \cos(\boldsymbol{k}_1 \cdot \boldsymbol{x}') + t_2 \cos(\boldsymbol{k}_2 \cdot \boldsymbol{x}') + \mathcal{O}(|\boldsymbol{t}|^2), \quad \boldsymbol{t} \in W.$$

*Furthermore, the solution depends analytically on $\boldsymbol{t} \in W$. In a neighbourhood of $(0, 0, \boldsymbol{c}^\star)$ in $((C^{1,\gamma}_{per,e}(\overline{\Omega^0}))^2 \times C^{1,\gamma}_{per,o}(\overline{\Omega^0})) \times C^{2,\gamma}_{per,e}(\mathbb{R}^2) \times \mathbb{R}^2$, these are the only non-trivial solutions, except for two two-parameter families of 2½-dimensional solutions which can be obtained by simple bifurcation from nearby points where the kernel of the linearisation is one-dimensional.*

**Remark 4.2.** Propositions 3.1 and 3.3 show that it is possible to satisfy the assumptions of the theorem. Indeed for any $\alpha > 0$ and $d > 0$, we can choose the lengths $|\boldsymbol{k}_1|$ and $|\boldsymbol{k}_2|$ so that $\sqrt{\alpha^2 - n^2|\boldsymbol{k}_j|^2} \notin \frac{\pi}{d}\mathbb{Z}_+$ and $\kappa(|\boldsymbol{k}_j|) > 0$ for $j = 1, 2$ and $n \in \mathbb{Z}$. Since $\gamma_j$ depends only on the lengths $|\boldsymbol{k}_j|$, we can then choose $\theta$, the angle between $\boldsymbol{k}_1$ and $\boldsymbol{k}_2$, so that (3.8) is fulfilled. We also choose it so that (2.7) is satisfied, that is

$$2n_1n_2|\boldsymbol{k}_1||\boldsymbol{k}_2|\cos(\theta) \neq \alpha^2 - n_1^2|\boldsymbol{k}_1|^2 - n_2^2|\boldsymbol{k}_2|^2 - \frac{n_3^2\pi^2}{d^2}.$$

for all $n_1, n_2 \in \mathbb{Z}$ and $n_3 \in \mathbb{Z}_+$. Note that we only have to avoid finitely many angles and that the case when either $n_1$ or $n_2$ vanishes is already satisfied by the choice of $|\boldsymbol{k}_1|$ and $|\boldsymbol{k}_2|$. Then there exists a $\boldsymbol{c} = \boldsymbol{c}^\star$ such that $\boldsymbol{k}_1$ and $\boldsymbol{k}_2$ are roots

of the dispersion equation. Moreover, the transversality condition is satisfied. At least in the case of a symmetric lattice, $|\boldsymbol{k}_1| = |\boldsymbol{k}_2|$, we can then choose $\sigma$ outside some countable set, which depends only on $\boldsymbol{k}_1$, $\boldsymbol{k}_2$, $\alpha$ and $d$, such that the dispersion equation has exactly the four roots $\pm\boldsymbol{k}_1$, $\pm\boldsymbol{k}_2$ in the lattice $\Lambda'$.

In the irrotational case, $\alpha = 0$, the first and third conditions are always satisfied (see Remark 3.5). As above, the second condition can be verified outside an exceptional set of parameter values in the case of a symmetric lattice; see Remark 3.2 and Reeder & Shinbrot [31].

**Remark 4.3.** In the proof of this result we will work with the reduced equation (2.6) for the surface profile. We loose a little bit of generality in doing this since we have to impose the non-resonance condition (2.7) which might be a bit stronger than needed if we were to work directly with problem (2.5) (see Cases III and IV in Section 3.1). On the other hand, we gain the elegance of the reduced equation and the bifurcation conditions become simpler to state.

Before giving the proof of the theorem we start with some technical lemmas. In what follows it will be useful to introduce the notation

$$X^k := C^{k,\gamma}_{per,e}(\mathbb{R}^2).$$

We rewrite (2.6) in the form

$$H(\eta, \boldsymbol{c}) := \sum_{j=1}^{2} c_j [D_\eta v_j[0, \boldsymbol{c}](\eta)]_{z=0} + g\eta - \sigma\Delta\eta - S(\eta, \boldsymbol{c}) = 0,$$

where

$$S(\eta, \boldsymbol{c}) := R(\boldsymbol{v}(\eta, \boldsymbol{c}), \eta) - \sum_{j=1}^{2} c_j([v_j(\eta, \boldsymbol{c}) - D_\eta v_j[0, \boldsymbol{c}](\eta)]_{z=0})$$

satisfies $S(0, \boldsymbol{c}) = 0$ and $D_\eta S[0, \boldsymbol{c}] = 0$. Note that

$$H\colon B_{r_0}(0; X^2) \times \mathbb{R}^2 \to X^0$$

is analytic with respect to $\eta$ and $\boldsymbol{c}$ by Theorem 2.1. It is convenient to represent

$$\begin{aligned} H(\eta, \boldsymbol{c}) &= D_\eta H[0, \boldsymbol{c}^\star](\eta) + D^2_{\eta, c_1} H[0, \boldsymbol{c}^\star](\eta, c_1 - c_1^\star) \\ &\quad + D^2_{\eta, c_2} H[0, \boldsymbol{c}^\star](\eta, c_2 - c_2^\star) + H_r(\eta, \boldsymbol{c}) \\ &= L(\eta) + (c_1 - c_1^\star)L_1(\eta) + (c_2 - c_2^\star)L_2(\eta) + H_r(\eta, \boldsymbol{c}). \end{aligned}$$

Here $L$, $L_1$ and $L_2$ are linear operators of $\eta$, given by

$$L(\eta) = \sum_{j=1}^{2} c_j^\star [D_\eta v_j[0, \boldsymbol{c}^\star](\eta)]_{z=0} + g\eta - \sigma\Delta\eta$$

and

$$\begin{aligned}L_1(\eta) &= \partial_{c_1}\left(\sum_{j=1}^{2} c_j[D_\eta v_j[0,\boldsymbol{c}](\eta)]_{z=0}\right)\Bigg|_{\boldsymbol{c}=\boldsymbol{c}^\star}\\ &= [D_\eta v_1[0,\boldsymbol{c}^\star](\eta)]_{z=0} + \sum_{j=1}^{2} c_j^\star[D_\eta v_j[0,c_1-c_1^\star,0](\eta)]_{z=0},\\ L_2(\eta) &= \partial_{c_2}\left(\sum_{j=1}^{2} c_j[D_\eta v_j[0,\boldsymbol{c}](\eta)]_{z=0}\right)\Bigg|_{\boldsymbol{c}=\boldsymbol{c}^\star}\\ &= [D_\eta v_2[0,\boldsymbol{c}^\star](\eta)]_{z=0} + \sum_{j=1}^{2} c_j^\star[D_\eta v_j[0,0,c_2-c_2^\star](\eta)]_{z=0},\end{aligned}$$

(note that $D_\eta \boldsymbol{v}[0,\boldsymbol{c}](\eta)$ is linear in $\boldsymbol{c}$), while the remainder

$$H_r(\eta,\boldsymbol{c}) := \sum_{j=1}^{2}(c_j-c_j^\star)[D_\eta v_j[0,\boldsymbol{c}-\boldsymbol{c}^\star](\eta)]_{z=0} - S(\eta,\boldsymbol{c}) \tag{4.1}$$

satisfies the estimate

$$\|H_r(\eta,\boldsymbol{c})\|_{X^0} \le C(\|\eta\|_{X^2} + |\boldsymbol{c}-\boldsymbol{c}^\star|^2)\|\eta\|_{X^2}.$$

We will later use the fact that if $\eta$ is constant in the direction $\boldsymbol{\lambda}_j$, then so are $L(\eta)$, $L_1(\eta)$, $L_2(\eta)$, $H_r(\eta,\boldsymbol{c})$, and hence $H(\eta,\boldsymbol{c})$. Indeed, this follows directly from the definitions and the last part of Theorem 2.1.

Let us study the action of the operator $D_\eta H[0,\boldsymbol{c}]$ in terms of the Fourier coefficients of $\eta$. Abbreviating the Fourier coefficients of $D_\eta H[0,\boldsymbol{c}]$ to $\widehat{D_\eta H}^{(\boldsymbol{k})}$, we find that

$$\widehat{D_\eta H}^{(\boldsymbol{k})} = \rho(\boldsymbol{c},\boldsymbol{k})\hat{\eta}^{(\boldsymbol{k})}, \quad \boldsymbol{k}\in\Lambda'\setminus\{0\},$$

where

$$\rho(\boldsymbol{c},\boldsymbol{k}) = g + \sigma|\boldsymbol{k}|^2 - \frac{(\boldsymbol{c}\cdot\boldsymbol{k})^2}{|\boldsymbol{k}|^2}\phi'(0;|\boldsymbol{k}|) + \alpha\frac{(\boldsymbol{c}\cdot\boldsymbol{k})(\boldsymbol{c}\cdot\boldsymbol{k}_\perp)}{|\boldsymbol{k}|^2}$$

is the expression from the dispersion equation (3.3), and

$$\widehat{D_\eta H}^{(0)} = g\hat{\eta}^{(0)}.$$

Similarly, we find that

$$\widehat{L(\eta)}^{(\boldsymbol{k})} = \rho(\boldsymbol{c}^\star,\boldsymbol{k})\hat{\eta}^{(\boldsymbol{k})}, \quad \widehat{L_j(\eta)}^{(\boldsymbol{k})} = \partial_{c_j}\rho(\boldsymbol{c}^\star,\boldsymbol{k})\hat{\eta}^{(\boldsymbol{k})}, \quad \boldsymbol{k}\in\Lambda'\setminus\{0\}.$$

By the assumptions of the theorem, $\rho(\boldsymbol{c}^\star,\boldsymbol{k}_j)=0$, $j=1,2$ and $\rho(\boldsymbol{c}^\star,\boldsymbol{k})\neq 0$, $\boldsymbol{k}\neq\pm\boldsymbol{k}_1, \boldsymbol{k}\neq\pm\boldsymbol{k}_2$. Using this together with standard properties of Fourier multiplier operators on Hölder spaces (see e.g. Bahouri, Chemin & Danchin [3, Prop. 2.78]), one obtains the following lemma.

**Lemma 4.4.** *The operator* $L\colon X^2 \to X^0$ *is Fredholm of index 0. Its kernel* $\ker L$ *is two-dimensional and spanned by the functions*

$$\eta_j(\boldsymbol{x}') = \cos(\boldsymbol{k}_j \cdot \boldsymbol{x}'), \quad j = 1, 2,$$

*and the operator* $L\colon \tilde{X}^2 \to \tilde{X}^0$ *is invertible, where* $\tilde{X}^k$ *denotes the orthogonal complement of* $\ker L$ *in* $X^k$ *with respect to the* $L^2_{per}$ *inner product.*

This allows us to use the Lyapunov-Schmidt reduction. For this purpose we split

$$\eta = t_1\eta_1 + t_2\eta_2 + \tilde{\eta}, \tag{4.2}$$

where $\tilde{\eta} \in \tilde{X}^2$. Thus the problem is written as

$$H(t_1\eta_1 + t_2\eta_2 + \tilde{\eta}, \boldsymbol{c}) = 0. \tag{4.3}$$

Let $P_j$ be the orthogonal projection in $L^2_{per}(\mathbb{R}^2)$ on the one-dimensional subspace spanned by $\eta_j$ and let $\tilde{P} = I - \sum_{j=1}^2 P_j$. Taking projections in (4.3) we obtain the $2 \times 2$ system

$$\begin{aligned}
&(c_1 - c_1^\star)t_1 P_1 L_1(\eta_1) + (c_2 - c_2^\star)t_1 P_1 L_2(\eta_1) + P_1 H_r(t_1\eta_1 + t_2\eta_2 + \tilde{\eta}, \boldsymbol{c}) = 0,\\
&(c_1 - c_1^\star)t_2 P_2 L_1(\eta_2) + (c_2 - c_2^\star)t_2 P_2 L_2(\eta_2) + P_2 H_r(t_1\eta_1 + t_2\eta_2 + \tilde{\eta}, \boldsymbol{c}) = 0,
\end{aligned} \tag{4.4}$$

and an equation for the orthogonal part

$$L(\tilde{\eta}) + (c_1 - c_1^\star)L_1(\tilde{\eta}) + (c_2 - c_2^\star)L_2(\tilde{\eta}) + \tilde{P} H_r(t_1\eta_1 + t_2\eta_2 + \tilde{\eta}, \boldsymbol{c}) = 0. \tag{4.5}$$

Applying the implicit function theorem to (4.5), noting that $L$ is an isomorphism from $\tilde{X}^2$ to $\tilde{X}^0$, we obtain the following reduction.

**Lemma 4.5.** *There exist constants* $\tilde{\epsilon}, \tilde{\delta}_0 > 0$*, a neighbourhood* $V$ *of the origin in* $\tilde{X}^2$ *and a function* $\psi\colon B_{\tilde{\epsilon}}(0; \mathbb{R}^2) \times B_{\tilde{\delta}_0}(\boldsymbol{c}^\star; \mathbb{R}^2) \to V$ *such that* (4.2)*, with* $(\boldsymbol{t}, \boldsymbol{c}, \tilde{\eta}) \in B_{\tilde{\epsilon}}(0; \mathbb{R}^2) \times B_{\tilde{\delta}_0}(\boldsymbol{c}^\star; \mathbb{R}^2) \times V$ *solves* (4.3) *if and only if* $\tilde{\eta} = \psi(\boldsymbol{t}, \boldsymbol{c})$ *and* $t_1\eta_1 + t_2\eta_2 \in \ker L$ *solves the finite-dimensional problem* (4.4) *with* $\tilde{\eta} = \psi(\boldsymbol{t}, \boldsymbol{c})$*. The function* $\psi$ *has the properties* $\psi(0, \boldsymbol{c}) = 0$ *and* $D_{\boldsymbol{t}}\psi[0, \boldsymbol{c}] = 0$.

For convenience we write $\tilde{\eta}(\boldsymbol{t}, \boldsymbol{c})$ instead of $\psi(\boldsymbol{t}, \boldsymbol{c})$ below. Note that $\tilde{\eta}(0, t_2, \boldsymbol{c})$ is constant in the direction $\boldsymbol{\lambda}_1$ and therefore independent of $\boldsymbol{k}_1 \cdot \boldsymbol{x}'$. Indeed, this follows by repeating the analysis in subspaces consisting of functions which only depend on $\boldsymbol{k}_2 \cdot \boldsymbol{x}'$ and using the mapping properties of the operators involved. Similarly, $\tilde{\eta}(t_1, 0, \boldsymbol{c})$ is independent of $\boldsymbol{k}_2 \cdot \boldsymbol{x}'$. We need one more technical lemma before the proof of the main result.

**Lemma 4.6.** *The remainder term $H_r$ satisfies*

$$\begin{aligned} H_r(t_1\eta_1 + t_2\eta_2 + \tilde{\eta}(\boldsymbol{t}, \boldsymbol{c}), \boldsymbol{c}) = a_1 t_1^2 \cos(2\boldsymbol{k}_1 \cdot \boldsymbol{x}') + a_2 t_1 t_2 \cos((\boldsymbol{k}_1 + \boldsymbol{k}_2) \cdot \boldsymbol{x}') \\ + a_3 t_1 t_2 \cos((\boldsymbol{k}_1 - \boldsymbol{k}_2) \cdot \boldsymbol{x}') + a_4 t_2^2 \cos(2\boldsymbol{k}_2 \cdot \boldsymbol{x}') \\ + a_5 t_1^2 + a_6 t_2^2 + \mathcal{O}(|\boldsymbol{t}|(|\boldsymbol{t}|^2 + |\boldsymbol{c} - \boldsymbol{c}^\star|^2)) \end{aligned} \tag{4.6}$$

*as $(\boldsymbol{t}, \boldsymbol{c} - \boldsymbol{c}^\star) \to 0$, where $a_j = a_j(\boldsymbol{c})$, $j = 1, \ldots, 6$, are real constants depending on $\boldsymbol{c}$.*

**Proof.** Considering formula (4.1) for $H_r$, we see that the first term is quadratic in $\boldsymbol{c} - \boldsymbol{c}^\star$ and linear in $\eta$. Therefore, when replacing $\eta$ by $t_1\eta_1 + t_2\eta_2 + \tilde{\eta}(\boldsymbol{t}, \boldsymbol{c})$ we get a contribution to the remainder of order $\mathcal{O}(|\boldsymbol{t}||\boldsymbol{c} - \boldsymbol{c}^\star|^2)$. The second term in (4.1) is quadratic in $\eta$ to lowest order. The quadratic part is obtained by forming products of terms involving $e^{\pm i\boldsymbol{k}_j \cdot \boldsymbol{x}'}$ and by evenness and realness we obtain precisely the above expression. □

**Proof of Theorem 4.1.** We think of the reduced system (4.4) with $\tilde{\eta} = \tilde{\eta}(\boldsymbol{t}, \boldsymbol{c})$ as a system of two scalar equations with respect to $c_1$ and $c_2$, while $t_1, t_2$ are parameters. Abusing notation, we identify the projection $P_j f$ of $f \in X^0$ with the coefficient in front of $\eta_j$ in its cosine series. We then obtain a $2 \times 2$ system of scalar equations

$$\begin{aligned} t_1\nu_{11}(c_1 - c_1^\star) + t_1\nu_{12}(c_2 - c_2^\star) + P_1 H_r(t_1\eta_1 + t_2\eta_2 + \tilde{\eta}(\boldsymbol{t}, \boldsymbol{c}), \boldsymbol{c}) = 0, \\ t_2\nu_{21}(c_1 - c_1^\star) + t_2\nu_{22}(c_2 - c_2^\star) + P_2 H_r(t_1\eta_1 + t_2\eta_2 + \tilde{\eta}(\boldsymbol{t}, \boldsymbol{c}), \boldsymbol{c}) = 0, \end{aligned}$$

where

$$\nu_{lj} = P_l L_j(\eta_l) = \partial_{c_j} \rho(\boldsymbol{c}^\star, \boldsymbol{k}_l).$$

We next note that

$$P_j \left( \left[ H_r(t_1\eta_1 + t_2\eta_2 + \tilde{\eta}(\boldsymbol{t}, \boldsymbol{c}), \boldsymbol{c}) \right]_{t_j = 0} \right) = 0, \quad j = 1, 2. \tag{4.7}$$

Indeed, as remarked after Lemma 4.5, $\tilde{\eta}(0, t_2, \boldsymbol{c})$ depends only on $\boldsymbol{k}_2 \cdot \boldsymbol{x}'$ and the same is true for

$$H_r(t_2\eta_2 + \tilde{\eta}(0, t_2, \boldsymbol{c}), \boldsymbol{c}).$$

Thus, its cosine series does not contain the mode $\eta_1$, which explains (4.7). A similar argument works for $j = 2$. It follows that

$$P_j \left( H_r(t_1\eta_1 + t_2\eta_2 + \tilde{\eta}(\boldsymbol{t}, \boldsymbol{c}), \boldsymbol{c}) \right) = t_j \Psi_j(\boldsymbol{t}, \boldsymbol{c}), \quad j = 1, 2,$$

where $\Psi_j$ is analytic in $B_{\tilde{\epsilon}}(0; \mathbb{R}^2) \times B_{\tilde{\delta}_0}(\boldsymbol{c}^\star; \mathbb{R}^2)$ with $\Psi_j(\boldsymbol{t}, \boldsymbol{c}) = \mathcal{O}(|\boldsymbol{t}|^2 + |\boldsymbol{c} - \boldsymbol{c}^\star|^2)$ in view of (4.6). We can thus rewrite the system as

$$\begin{aligned} t_1(\nu_{11}(c_1 - c_1^\star) + \nu_{12}(c_2 - c_2^\star) + \Psi_1(\boldsymbol{t}, \boldsymbol{c})) = 0, \\ t_2(\nu_{21}(c_1 - c_1^\star) + \nu_{22}(c_2 - c_2^\star) + \Psi_2(\boldsymbol{t}, \boldsymbol{c})) = 0. \end{aligned} \tag{4.8}$$

We obtain a solution by solving the system

$$\nu_{11}(c_1 - c_1^\star) + \nu_{12}(c_2 - c_2^\star) = -\Psi_1(\boldsymbol{t}, \boldsymbol{c}),$$
$$\nu_{21}(c_1 - c_1^\star) + \nu_{22}(c_2 - c_2^\star) = -\Psi_2(\boldsymbol{t}, \boldsymbol{c}).$$

The determinant of the coefficients in the left-hand side is

$$\partial_{c_1}\rho(\boldsymbol{c}^\star, \boldsymbol{k}_1)\partial_{c_2}\rho(\boldsymbol{c}^\star, \boldsymbol{k}_2) - \partial_{c_1}\rho(\boldsymbol{c}^\star, \boldsymbol{k}_2)\partial_{c_2}\rho(\boldsymbol{c}^\star, \boldsymbol{k}_1).$$

This is non-zero if and only if the transversality condition (3.12) is fulfilled. Since the right-hand side is $\mathcal{O}(|\boldsymbol{t}|^2 + |\boldsymbol{c} - \boldsymbol{c}^\star|^2)$, it follows from the implicit function theorem that we can solve this system for $\boldsymbol{c}$ if $|\boldsymbol{t}| < \epsilon$ for some $\epsilon > 0$, and that

$$\boldsymbol{c}(\boldsymbol{t}) = \boldsymbol{c}^\star + \mathcal{O}(|\boldsymbol{t}|^2).$$

This results in a two-parameter family of solutions

$$\eta(\boldsymbol{x}'; \boldsymbol{t}) = t_1\eta_1(\boldsymbol{x}') + t_2\eta_2(\boldsymbol{x}') + \tilde{\eta}(\boldsymbol{x}'; \boldsymbol{t}, \boldsymbol{c}(\boldsymbol{t})), \quad |\boldsymbol{t}| < \epsilon,$$

for which the vector field

$$\boldsymbol{v}(\boldsymbol{x}; \eta(\boldsymbol{t}), \boldsymbol{c}(\boldsymbol{t}))$$

solves problem (2.5).

Let us discuss the special case when $t_1 = 0$. If so, then (4.8) reduces to the scalar equation

$$\nu_{21}(c_1 - c_1^\star) + \nu_{22}(c_2 - c_2^\star) = -\Psi_2(0, t_2, \boldsymbol{c})$$

(and the additional trivial solution $t_1 = t_2 = 0$). Under the non-degeneracy condition at least one of the coefficients $\nu_{2j}$ in the left-hand side is non-zero, and we may solve for the corresponding parameter $c_j$, thus obtaining a two-dimensional family of non-trivial solutions parametrised by $t_2$ and the other parameter $c_{j'}$. We claim that these solutions are 2½-dimensional. Indeed, this follows since we could repeat the analysis in a space of functions that only depend on $\bar{x} = \boldsymbol{x}' \cdot \boldsymbol{k}_j$. The same of course applies if $t_2 = 0$. These solutions can also be obtained by one-dimensional bifurcations from points $\boldsymbol{c}$ close but not equal to $\boldsymbol{c}^\star$, where the kernel is spanned by either $\eta_1$ or $\eta_2$. The two-parameter family of genuinely three-dimensional solutions parametrised by $t_1$ and $t_2$ intersects these two-parameter families of 2½-dimensional waves along the curves $t_1 = 0$, $\boldsymbol{c} = \boldsymbol{c}(0, t_2)$ and $t_2 = 0$, $\boldsymbol{c} = \boldsymbol{c}(t_1, 0)$ in the four dimensional parameter space $(t_1, t_2, c_1, c_2) \in \mathbb{R}^4$. Thus the family of genuinely three-dimensional waves connects different 2½-dimensional waves through what can be described as 'dimension breaking bifurcations' (see e.g. GROVES, HARAGUS & SUN [19] and GROVES, SUN & WAHLÉN [22]). This completes the proof. □

**Remark 4.7.** The last observation in the proof holds more generally. If we have a kernel spanned by $\cos(\boldsymbol{x}' \cdot \boldsymbol{k}_0)$, $\boldsymbol{k}_0 \in \Lambda'$, for some value of $\boldsymbol{c}$, then we can fix one of the parameters and apply local bifurcation theory with a one-dimensional kernel. The resulting family of solutions will be 2½-dimensional since we can consider the same problem but for functions that depend only on $\bar{x} = \boldsymbol{x}' \cdot \boldsymbol{k}_0$. This explains why one has to consider at least two-dimensional kernels in order to construct genuinely three-dimensional waves.

*Acknowledgements.* Open access funding provided by Lund University. This project has received funding from the European Research Council (ERC) under the European Union's Horizon 2020 research and innovation programme (grant agreement no 678698). The authors are grateful to the referees for constructive criticism.

Doubly Periodic Water Waves with Vorticity

E. Lokharu
Department of Mathematics,
Linköping University,
581 83 Linköping
Sweden.

and

D. S. Seth & E. Wahlén
Centre for Mathematical Sciences,
Lund University,
PO Box 118, 22100 Lund
Sweden.
e-mail: erik.wahlen@math.lu.se